\documentclass{amsproc}
\usepackage{amssymb,amscd}

\def\tensor{\otimes}
\def\isomorphic{\cong}
\DeclareMathOperator{\codim}{codim}

\DeclareMathOperator{\Sym}{Sym}

\def\directsum{\oplus}
\def\N{\operatorname{N}}

\def\union{\cup}

\def\intersect{\cap}

\DeclareMathOperator{\Gal}{Gal}
\def\Pr{{\bf P}} %projective space
\def\k{k} %base field
\def\kbar{{\overline{\k}}} %its algebraic closure
\def\X{X} %variety
\def\C{C} %curve
\def\OC{\mathcal{O}_{\C}} %structure sheaf
\def\SymdC{{\Sym^\d(\C)}} % dth symmetric power of C
\def\L{\mathcal{L}} %ample line bundle
\def\bundleE{\mathcal{E}} %vector bundle
\def\G{{\bf G}} %Grassmannian
\def\Gd{\G_\d} %Grassmannian parametrizing codim d subspaces
\def\V{V} %vecspace
\def\Vstar{\V^*} %its dual
\def\x{x} %point on X
\def\s{s} %section of L
\def\t{t} %another section
\def\f{f} %a hyperplane section for the main flipping algorithm.
\def\g{g} %genus
\def\N{N} %degree of line bundle L on C
\def\D{D} %divisor
\def\E{E} %another divisor
\def\d{d} %codimension d subspaces, and degree of D
\def\e{e} %degree of E
\def\HoOC#1{H^0\bigl(\OC(#1)\bigr)}
\def\HoL{H^0(\L)} %H^0(C,L)
\def\HoLone{H^0(\L_1)} %H^0(C,L_1)
\def\HoLtwo{H^0(\L_2)} %H^0(C,L_2)
\def\P{P} %point on C
\def\mP{m_\P} %order of vanishing at P
\def\W{W} %subspace in Grassmannian
\def\WD{\W_\D} %image of the divisor D
\def\WDstar{\W_\D^*}
\def\WE{\W_\E} %image of the divisor E
\def\WEstar{\W_\E^*}
\def\Wstar{\W^*}
\def\HoLm#1{H^0(\L - #1)} % \Holm(\D) is \WD
\def\HoLmD{{\HoLm{\D}}}
\def\J{J} % The Jacobian variety of C
\def\mul{\mu} %multiplication of sections
\def\HotwoL{H^0(2\L)} %square of the line bundle L
\def\HotwoLm#1{H^0(2\L - #1)} %part of the above.
\def\A{\mathcal{A}} %projective coordinate ring
\def\I{\mathcal{I}} %ideal in A
\def\Hyp{H} %hyperplane (f=0)
\def\linequiv{\sim}
\def\thedim{\delta} %dimensions of vector spaces being used
\def\x{x} %point on Jacobian
\def\xD{\x_\D}
\def\xE{\x_\E}

\newtheorem{theorem}{Theorem}[section]
\newtheorem{lemma}[theorem]{Lemma}

\newtheorem{lemmalg}[theorem]{Lemma/Algorithm}

\newtheorem{propalg}[theorem]{Proposition/Algorithm}
\newtheorem{thmalg}[theorem]{Theorem/Algorithm}

\theoremstyle{definition}
\newtheorem{definition}[theorem]{Definition}

\theoremstyle{remark}
\newtheorem{remark}[theorem]{Remark}

\numberwithin{equation}{section}

\begin{document}

\title{Linear algebra algorithms for divisors on an algebraic curve}

% author one information
\author{Kamal Khuri-Makdisi}
\address{Mathematics Department and Center for Advanced Mathematical Sciences,
American University of Beirut, Bliss Street, Beirut, Lebanon}
\email{kmakdisi@aub.edu.lb}

% author two information
%\author{}
%\address{}
%\email{}
%\thanks{}

\thanks{Version submitted for publication, March 26, 2002}

\begin{abstract}
We use an embedding of the symmetric $d$th power of any algebraic
curve $C$ of genus $g$ into a Grassmannian space to give algorithms
for working with divisors on $C$, using only linear algebra in vector
spaces of dimension $O(g)$, and matrices of size $O(g^2)\times O(g)$.
When the base field $k$ is finite, or if $C$ has a rational point over
$k$, these give algorithms for working on the Jacobian of $C$ that
require $O(g^4)$ field operations, arising from the Gaussian
elimination.  Our point of view is strongly geometric, and our
representation of points on the Jacobian is fairly simple to deal
with; in particular, none of our algorithms involves arithmetic with
polynomials.  We note that our algorithms have the same asymptotic
complexity for general curves as the more algebraic algorithms in
Hess' 1999 Ph{.}D{.} thesis \cite{HessThesis}, which works with
function fields as extensions of $k[x]$.  However, for special classes
of curves, Hess' algorithms are asymptotically more efficient than
ours, generalizing other known efficient algorithms for special
classes of curves, such as hyperelliptic curves~\cite{Cantor},
superelliptic curves~\cite{GPS}, and $C_{ab}$ curves~\cite{Suzuki}; in
all those cases, one can attain a complexity of $O(g^2)$.
\end{abstract}

\maketitle

\section{Introduction}
\label{section1}

Let $\C$ be a smooth algebraic curve of genus $\g$ over a field $\k$, and
assume for simplicity in this introduction that $\C$ has a rational point
over $\k$.  When $\g = 1$, $\C$ is an elliptic curve, and can be
represented as a plane cubic; this is the Weierstrass model.  In that case,
the group law on $\C$ is easy to describe and to implement, and has led to
many effective algorithms in cryptography that make use of the difficulty
of solving the discrete logarithm problem in the group of $\k$-rational
points of $\C$ when $\k$ is a finite field.  If the genus of $\C$ is
larger, then one can try to work with the group law on the Jacobian of $\C$.
If $\C$ admits a morphism of small degree to $\Pr^1$ (for instance, if $\C$
is hyperelliptic), and if the genus $\g$ is large compared
to the order of $\k$, then subexponential algorithms are known for the
discrete logarithm on the Jacobian; see \cite{AdDMHu} for hyperelliptic
curves, \cite{GPS} for superelliptic curves, and \cite{HessThesis} for
general curves.  Computing in the Jacobian of a large genus curve is
nonetheless of intrinsic interest.

One problem is that for large $\g$, the Jacobian is somewhat difficult to
describe directly as an algebraic variety.  For instance, if one wishes to
embed the Jacobian into projective space, the standard way is to take a
very ample line bundle (say, $3$ or $4$ times the theta-divisor), which
leads to an embedding of the Jacobian into a projective space of
exponentially large dimension (e.g., $3^\g - 1$ or $4^\g - 1$).  This is
not practical for explicit computations, especially as the equations
defining the group law on the Jacobian will probably be similarly
intractable.
% (A different way to embed the Jacobian into a projective space of
% unmanageably large dimension is described in Chow [Reference?\qqq].)
The other way to deal with the Jacobian is to use the fact that it is
essentially an ideal class group attached to the function field of $\C$; a
divisor is a certain kind of ideal, and linear equivalence between divisors
is equivalence in the ideal class group.  Thus one works directly with
divisors (or ideals) on $\C$, and spends time mostly on finding appropriate
elements of the function field in order to reduce divisors to a canonical
form modulo linear equivalence.  If the curve $\C$ is hyperelliptic, this
is relatively easy, since the Jacobian is analogous to the ideal class
group of a quadratic field; questions about the Jacobian can be translated
into questions about binary quadratic forms over the polynomial ring
$\k[x]$.  This allows one to implement the group operations on the Jacobian
of a hyperelliptic curve in ``time'' proportional to $\g^2$, as
in~\cite{Cantor}.  By this we mean that the algorithms require $O(\g^2)$
operations in the field $\k$.  The articles \cite{GPS} and \cite{Suzuki}
give $O(\g^2)$ algorithms for a wider class of curves, the superelliptic
and $C_{ab}$ curves, using more sophisticated techniques to work with the
ideal classes; the implied constants are however quite large, and there is
some possibility that they could grow with the genus $\g$.  As for a
general curve of genus $\g$, there are algorithms by F. Hess (see
\cite{HessThesis} and \cite{HessArticle}), which both implement the group
operations in the Jacobian, and determine its group structure over a finite
field.  Hess' algorithms begin with a map of degree $n$ from the curve $\C$
to $\Pr^1$, and work with ideals in the function field of $\C$, viewed as a
degree $n$ extension of $k[x]$.  If $n$ is fixed or bounded, then Hess'
algorithms require $O(\g^2)$ field operations to implement the group law on
the Jacobian.  This generalizes the previously mentioned results on
hyperelliptic, superelliptic, and $C_{ab}$ curves.  However, the minimum
possible $n$ for a general curve of genus $\g$ is approximately $\g/2$ (see
\cite{GrifHar}, p.~261).  In this case (and, more generally, for $n =
O(\g)$), the complexity of Hess' algorithms rises to $O(\g^4)$ field
operations per group operation in the Jacobian.  Besides Hess' algorithms
for general curves, we note two earlier geometric algorithms that are due
to Huang and Ierardi~\cite{HuangIerardi} and Volcheck~\cite{Volcheck}.
Both algorithms rely on a description of $\C$ as a plane curve with
singularities, and the steps can become fairly involved (approximately
$O(g^7)$ field operations are needed per group operation in the Jacobian).

In this paper, we give a set of straightforward geometric algorithms
for the Jacobian of $\C$ that
involve only $O(\g^4)$ field operations.  Moreover, our algorithms are
quite simple to implement, as they only involve linear algebra (row
reduction) on matrices of size $O(\g^2) \times O(\g)$, and in vector spaces
of dimension $O(\g)$.  In particular, no polynomial algebra is ever used
explicitly ({\em a fortiori}, no Gr\"obner bases or field extensions
either).  One can however argue that the representation of $\C$ that we
use implicitly works with quadratic polynomials in $O(\g)$ variables.
Specifically, we take a line bundle $\L$ on $\C$, whose degree $\N$ is
approximately $6\g$, and work with the projective representation of $\C$
arising from $H^0(\C,\L)$.  In our setting, the homogeneous ideal defining
$\C$ is generated by quadratic polynomials (elements of $\Sym^2 H^0(\C,
\L)$), which can be thought of as the kernel of the multiplication map
\begin{equation}
\label{equation1.1}
\mul: H^0(\C, \L) \tensor H^0(\C, \L) \to \Sym^2 H^0(\C,\L) \to H^0(\C,
\L^{\tensor 2}).
\end{equation}
We note that the dimensions of the vector spaces $H^0(\C,\L)$ and $H^0(\C,
\L^{\tensor 2})$ are both $O(\g)$, being approximately $5g$ and $11g$,
respectively.  Then $\mul$ can be
described in terms of any suitable bases for these vector spaces.  This
article does not present any new ideas on how to efficiently construct
curves $\C$ represented in the above form, and, what is more, how to
construct curves along with some nontrivial elements of their Jacobians.
Nonetheless, the first question is easy to answer in principle: for
example, choose a polynomial $f(x,y)$ which describes a plane curve
with singularities that is a projection of $\C$, and which passes through a 
given nonsingular point $\P_0$ such as $(0,0)$; then we can let $\L =
\OC(\N\P_0)$, and use the algorithms of \cite{HessThesis} and of
\cite{HessArticle} (or perhaps those of \cite{HuangIerardi} or
of~\cite{Volcheck}) to compute the global sections of $\L$ and of
$\L^{\tensor 2}$, and hence the map $\mul$.  This must only be done once,
in obtaining the initial description of $\C$ in the form suitable for our
algorithms.  As for the second question above, we have not come up with any 
interesting methods to produce divisors on our curve $\C$, corresponding
to nontrivial points on the Jacobian.  The best that we can do is to
suggest choosing the plane curve to pass through a small number of given
points $\P_1, \dots, \P_n$, with $n$ small compared to $\g$.  Then we can
obtain divisors generated by the $\P_1, \dots \P_n$, but these divisors may 
not be ``typical,'' for all we know.  Another topic that we have not
addressed (but which is covered in \cite{HessThesis}) is the question
of efficiently computing the order, or, even better, the group
structure, of the Jacobian over a finite field.

The basic forms of our algorithms, given in Section~\ref{section4}, are
fairly easy to describe and to understand.  The algorithms involve only
subspaces of the two vector spaces $H^0(\C,\L)$ and $H^0(\C, \L^{\tensor
2})$ and the map $\mul$, and require one to solve various systems of linear 
equations in this setting.  In Section~\ref{section5}, and to a lesser
extent in Section~\ref{section4}, we give some faster, but slightly more
complicated, algorithms for working with divisors on $\C$ and with the
Jacobian, including an algorithm that essentially solves the general
Riemann-Roch problem for a divisor (written as a difference of effective
divisors) $\D_1 - \D_2$ in time $O(\max(\g, \deg \D_1, \deg \D_2)^4)$.
These faster algorithms still require time $O(\g^4)$, but work
with vector spaces of smaller dimension; for example, we can get away with
taking the line bundle $\L$ to have degree approximately $3\g$, but then
need to work (implicitly) with polynomials of higher degree in the
projective embedding of $\C$.  Nonetheless, the structure of our faster
algorithms is virtually identical to that of the basic forms of our
algorithms.  Our results follow from three main insights:
\begin{enumerate}
\item
We use a simple representation for effective divisors $\D$ on $\C$, in
terms of the projective embedding given by the line bundle $\L$ above.
Namely, so long as the degree $\d$ of $\D$ is small compared to that of
$\L$, we can identify $\D$ uniquely by the linear subspace $\WDstar$ of
projective space that is spanned by $\D$ (in a way that makes sense even if
$\D$ 
has points of multiplicity greater than~$1$).  Dually, we represent
$\D$ by the space $\WD \subset H^0(\C,\L)$ consisting of global sections
$\s \in H^0(\C,\L)$ that vanish at $\D$.  Since the effective divisor $\D$
corresponds to a point on the symmetric power $\Sym^{\d} \C$, we are
essentially working with an embedding of $\C$ into the Grassmannian
parametrizing codimension $\d$ subspaces of $H^0(\C,\L)$.  As discussed
in Section~\ref{section2} below, this is much better than passing to a
projective model of $\Sym^{\d} \C$ using the Pl\"ucker embedding.
Perhaps other computational issues in algebraic geometry can become simpler
if one is willing to work more extensively with Grassmannians and linear
algebra, instead of always sticking to projective space and polynomials.
The author hopes that this article will promote further investigation
into the topic.
\item
We no longer insist on reducing divisors to a canonical form, but instead
represent an element of the Jacobian as an effective divisor of a certain
degree $\d \geq g+1$.  The basic form of our algorithms in
Section~\ref{section4} can use any $\d \geq 2\g+1$, and an improvement in
Section~\ref{section5} uses $\d \geq \g+1$.  It is quite doable to bring
this down to $\d = 
\g$, but at the expense of some complications in the algorithms.  As a
result, a given element of the Jacobian can be represented by many
different subspaces $\WD$, where only the linear equivalence class of $\D$
is well-defined.  We of course give an algorithm to determine when two
divisors represent the same point in the Jacobian; it also takes time
$O(\g^4)$.  (More precisely, our algorithms take time $O(\max(\d,\g)^4)$, 
in case someone should wish to work with divisors of extremely large
degree.)
\item
We work in terms of multiplication of sections of line bundles, as well as
a form of division (Lemma~\ref{lemma2.1} and
Lemma/Algorithm~\ref{lemmalg2.2}).  This amounts to working directly in the
projective coordinate ring of $\C$, instead of the full polynomial algebra
of the projective space in which $\C$ lies.  This produces a substantial
time savings, since for fixed $\ell << \g$, the space of degree $\ell$
polynomials has dimension $O({\dim V}^\ell) = O(\g^\ell)$, while the
restriction of such polynomials to $\C$ is the space $H^0(\C, \L^{\tensor
\ell})$, whose dimension is $O(\ell \g)$.  Moreover, the multiplication and
division operations are the basis of our addition and flipping algorithms
(Theorems/Algorithms \ref{thmalg3.4} and~\ref{thmalg3.6}), which are
fundamental to our work.  The flipping algorithm is especially important:
given a divisor $\D$ on $\C$, and a hyperplane containing $\D$, it allows
us to find the ``complementary'' divisor $\D'$ such that $\D + \D'$ is the
intersection of $\C$ with the hyperplane in question.  The algorithm is
inspired by geometric considerations, as explained in the discussion
preceding Theorem/Algorithm~\ref{thmalg3.6}.  The flipping algorithm is
what allows us to pass to divisor classes; for example, flipping twice with
respect to different hyperplanes replaces $\D$ with another divisor
linearly equivalent to it.
\end{enumerate}

We conclude this introduction by pointing out the connection between our
way of representing divisors, and the way in Cantor's algorithm for
hyperelliptic curves; it should be possible to say something similar for
superelliptic and $C_{ab}$ curves.  For simplicity, we restrict to the case
where $\k$ does not have characteristic~$2$, and so consider a curve $\C :
y^2 = f(x)$, where $f(x) \in \k[x]$ has degree $2\g + 1$.  Then a reduced
(effective) divisor $\D$ of degree $\d \leq \g$ is represented in Cantor's
algorithm by the pair of polynomials $(a(x), b(x))$, with $a(x)$ monic of
degree $\d$, $b(x)$ of degree at most $\d - 1$, and $b(x)^2 \equiv
f(x) \pmod{a(x)}$.  If we factor $a(x) = \prod_i(x-x_i)$, possibly with
some repeated $x_i$, then $\D$ is the sum of the points $(x_i, b(x_i))$ on
the curve.  On the other hand, in our representation, we choose a line
bundle $\L$ and represent $\D$ by the global sections of $\L$ vanishing at
$\D$.  We can take $\L = \OC(\N\P_\infty)$, for $\N \geq 3\g + 1$, where
$\P_\infty$ is the ``point at infinity'' with $v_{\P_\infty} (x) = 2$ and
$v_{\P_\infty} (y) = 2\g + 1$.  To simplify, we shall take $\N = 2m$ to
be even.  Then a basis for $H^0(\C, \L)$, arranged in increasing order of
pole at $\P_\infty$, is $\{1, x, x^2, \dots, x^\g, y, x^{\g+1}, xy,
x^{\g+2}, x^2 y, \dots, x^m\}$.  Those sections that vanish at $\D$ are
then of the form $p(x) a(x) + q(x) (y - b(x))$, where $\deg p(x) \leq m -
d$ and $\deg q(x) \leq m - 1 - \g$.  Our working with these sections using
linear algebra can be thought of as performing polynomial operations on
$a(x)$ and $b(x)$, but with a bound on the degree of intermediate results.
(Of course, our algorithms also work in a more general setting.)  An
alternative way to deal with the above is to consider the affine curve $\C
- \P_0$; its affine coordinate ring is $\k[x,y]/(y^2 - f(x)) \isomorphic
\k[x] \directsum \k[x]y$, and the elements vanishing at $\D$ form an ideal
generated (even as a $\k[x]$ module) by the two elements $a(x)$ and
$y-b(x)$.

\section*{Acknowledgements}
I am grateful to K. Murty for posing to me the question of whether one
could work algorithmically with abelian varieties without using explicit
equations.  I would also like to thank M. Thaddeus for very useful
conversations on linear systems, in particular, for his referring me to
the article \cite{Lazar}, and for his suggestion to think
primarily in terms of the space $\WD$ and sections of $\L$, rather than its
dual $\WDstar$ and the projective embedding.  
My thanks go to C. Gasbarri for helpful comments and references.
I also gratefully acknowledge
the support of the Clay Mathematics Institute, which partially supported
this research by funding me as a Clay Scholar in the summer of 2000, and
which supported my participation in the Clay Mathematics Institute
Introductory Workshop in Algorithmic Number Theory at MSRI in August 2000.
Finally, I would like to thank the mathematics department at Columbia
University for its hospitality at the time of the final revisions to this
manuscript.

%\section*{Notation} \qqq
%%%%%%%%%%%%%%%%%%%%%%%%%%%%%%%%%%%%%%%%%%%%%%%%%%%%%%%%%%%%%%%%%%%%%%%%%%%%
\section{Embedding the symmetric power of a curve into a Grassmannian}
\label{section2}
Let $\X$ be an algebraic variety over a field $\k$; for the moment,
we take $\k$ to be algebraically closed.  It is standard to try to embed
$\X$ into projective space by choosing a very ample line bundle $\L$
on $\X$, and using $\L$ to get an embedding of $\X$ into
$\Pr(H^0(\X,\L))$.  Here $\Pr(\V)$ is the projective space that
parametrizes {\em codimension} one subspaces of $\V$.  The standard
embedding associates $\x \in \X$ to the subspace $\{ \s \in H^0(\X,\L) \mid 
\s$ vanishes at $\x \}$.  More generally, one can try to embed $\X$ into a
Grassmannian variety, using a suitably positive rank $\d$ vector bundle
$\bundleE$.  We then associate to $\x \in \X$ the analogous subspace
$\W_\x$ of global sections of $\bundleE$ vanishing at $\x$.  Under suitable
hypotheses, this gives us an embedding $\varphi: \X \to \Gd(\V)$, where
$\Gd(\V)$ parametrizes codimension $\d$ subspaces $\W$ of $\V$.  The map
$\varphi$ can be composed with the Pl\"ucker embedding $\Gd(\V) \to
\Pr(\wedge^\d \V)$ to obtain an embedding of $\X$ into projective space.
For practical purposes, this may be unhelpful, as the dimension of
$\wedge^\d \V$ is $\binom{\dim\V}{\d}$-dimensional, which is significantly
greater than $\dim \V$ once $\d \geq 2$.  If moreover $\d$ increases with
$\dim V$, the growth in the dimension can be exponential; for example,
$\binom{2n}{n}$ grows roughly as $4^n/\sqrt{\pi n}$ by Stirling's formula.
This suggests that we would be better off working with subspaces of $V$,
and doing explicit calculations using linear algebra, instead of passing to
a possibly exponential number of explicit projective coordinates via the
Pl\"ucker embedding.  The case that interests us most is when $\X$ is the
$\d$th symmetric power $\SymdC$ of a curve $\C$.  We shall embed $\SymdC$
into a certain Grassmannian without explicitly mentioning vector bundles.
The alert reader will see that we are essentially working with the $\bf
Quot$ scheme by hand (see Section 8.2 of \cite{BLR}).

Let $\g$ be the genus of $\C$, take $\d \geq 1$, and fix once and for all a
line bundle $\L$ on $\C$ of degree $\N \geq 2\g + \d \geq 2\g + 1$.  The
Riemann-Roch theorem implies that the space $\V = \HoL = H^0(\C,\L)$ has
dimension $\N + 1 - \g$, and $\L$ gives an embedding of $\C$ into $\Pr(\V)
\isomorphic \Pr^{\N-\g}$.  In fact, a well-known theorem by Castelnuovo,
Mattuck, and Mumford states that this embedding of $\C$ is projectively
normal, and there is an even stronger theorem due to Fujita and
Saint-Donat, which states that if in fact $\N \geq 2\g + 2$, then the ideal
of $\C$ is generated by quadrics (\cite{Lazar}, Section 1.1).  Now a
point on $\SymdC$ can be identified with an effective divisor $\D =
\sum_{\P}\mP \P$ on $\C$, of degree $\d = \sum \mP$; we associate to $\D$
the subspace $\WD$ given by
\begin{equation}
\label{equation2.1}
\WD = \HoLmD \subset \V = \HoL.
\end{equation}
Thus $\WD$ consists
of the global sections $\s$ of $\L$ vanishing at $\D$ (i.e.,
$\s$ vanishes at each $\P$ to order at least $\mP$).  Note that we are
freely abusing notation by often failing to distinguish between line
bundles and divisors, and by referring to tensor operations on line bundles
additively when it suits us; it would be more standard to write $\WD =
H^0(\L(-\D))$.  Due to the condition on the degree $\N$, the line bundle
$\L - \D$ is nonspecial and has no base points, so $\WD$ has codimension
$\d$ in $\V$, and we can recover $\D$ from $\WD$.  We summarize this
discussion in the following lemma, while noting incidentally that the
argument also extends to the trivial case $\d = 0$.

\begin{lemma}[Representation of divisors]
\label{lemma2.0.5}
Let $\C$ be an algebraic curve of genus $\g$, and let $\L$ be a line bundle
on $\C$ of degree $\N$.  Let $\D = \sum \mP \P$ be an effective divisor on
$\C$ of degree $\d$ (we allow $\D$ to be empty, in which case $\d = 0$).
Let $\WD = \HoLmD \subset \HoL = \V$.  Then if $\d \leq \N - 2\g + 1$, it
follows that $\WD$ has codimension $\d$ in $\V$.  If furthermore $\d \leq
\N - 2\g$, then $\HoLmD$ has no base points, in the sense that
$\HoLm{\D-\P} \subsetneq \HoLm{\D}$ for all points $\P$ of $\C$.  In that
case, the multiplicity $\mP$ can be recovered by noting that $\mP$ is the
minimum order to which any section $\s \in \WD$ vanishes at $\P$ (viewing
$\s$ as a global section of $\L$).
\end{lemma}

\begin{proof}
This follows immediately from the Riemann-Roch theorem.
\end{proof}

Thus the map $\D \mapsto \WD$ is injective, and allows us to view $\SymdC$
as a subset of $\Gd(\V)$.  For $\N$ large compared to $\d$ and $\g$,
Theorem~$8.2.8'$ of \cite{BLR} guarantees that this injection of
points is an embedding of varieties; we suspect that this holds in our
setting, but have not checked the matter further, as we do not need
this fact.  Another way to  
view this map is via the embedding of $\C$ into $\Pr(\V)$ given by our
choice of $\L$.  Given a degree $\d$ divisor $\D$ on $\C$, then its points
span a certain $\d$-dimensional subspace $\WDstar$ of the dual space
$\Vstar$ (now identifying points of $\Pr(\V)$ with nonzero elements of
$\Vstar$ up to multiplication of scalars).  Assuming that the
multiplicities in $\D$ are all at most 1 (i.e., that $\D = \sum_i \P_i$ for
distinct $\P_i$ can be viewed as a reduced scheme), then $\WDstar$ is
precisely the annihilator of $\WD$ with respect to the pairing between $\V$
and $\Vstar$.  If $\D$ has some higher multiplicities, we can still make
sense of this last statement by understanding $\WDstar$ to be spanned not
only by $\P$, but by the $\mP$th infinitesimal neighborhood (i.e.,
osculating plane) of $\P$ on $\C$ (see \cite{GrifHar}, pp.~248 and~264).
We shall interchangeably use $\WDstar$ to refer to the $\d$-dimensional
subspace of $\Vstar$ and to the corresponding $(\d-1)$-plane in
$\Pr(\V)$.  In this point of view, the fact that $\HoLmD$ has no base
points means that the points of $\D$ are linearly independent (even
including multiplicity), and that the plane $\WDstar$ intersects $\C$
precisely in the divisor $\D$.  This gives another way to see the
injectivity of the map $\D \mapsto \WD$.

So far, we have worked over an algebraically closed field.  The above
discussion continues to hold if $\C$ and $\L$ are defined over an arbitrary
perfect base field $\k$, provided that we restrict ourselves to divisors
$\D$ rational over $\k$. This means that for all Galois automorphisms
$\sigma \in \Gal(\kbar/\k)$, we have an equality of divisors (and
not just of divisor classes) $\sigma \D = \D$.  In that case, the vector
spaces $\V$ and $\WD$ are finite-dimensional over $\k$, and all our work
can be done over the field $\k$.  This is quite convenient computationally, 
as we do not need to work with individual points of $D$, which can
themselves be defined over extensions of $\k$.  Still, our main goal in
working with $\SymdC$ is to descend to the Jacobian $\J$ of $\C$.  In that
setting, the above notion of rationality can fall somewhat short: a point
of $\J(\k)$ need not be representable by a divisor rational over $\k$,
since we now only require equality of the divisor {\em classes} $\D$ and
$\sigma(\D)$.  This problem arises only if $\C$ itself has no points
defined over $\k$; in that case, the obstruction to lifting a point of
$\J(\k)$ to a divisor rational over $\k$ lies in the Brauer group of $\k$
(see Remark~1.6 of \cite{Milne}).  We are primarily interested in the case
where $\k$ is a finite field (so the Brauer group is trivial), so we shall
henceforth ignore the
above distinction between rationality questions on the symmetric product
and on the Jacobian.  We merely remark that in working with various spaces
of sections of line bundles (such as $\HoLmD$), we shall work over the
algebraic closure $\kbar$ when needed, and then deduce the desired results
for $\k$.  A typical example is the following lemma, which is a slight
generalization of the theorem by Castelnuovo, Mattuck, and Mumford cited
above.

\begin{lemma}[``Multiplication'']
\label{lemma2.1}
Let $\C$ be a smooth, geometrically connected algebraic curve of genus $\g$
over a perfect field $\k$, and let $\L_1$ and $\L_2$ be two line bundles on
$\C$ (both defined over $\k$), of 
degrees $\d_1$ and $\d_2$.  Assume that $\d_1, d_2 \geq 2\g + 1$.
Then the canonical (multiplication) map
\begin{equation}
\label{equation2.2}
\mul : \HoLone \tensor_{\k} \HoLtwo \to H^0(\L_1 \tensor \L_2)
\end{equation}
is surjective.
\end{lemma}

\begin{proof}
We sketch a proof, adapting the exposition in Sections 1.3 and~1.4 of
\cite{Lazar}, which treats the case $\L_1 = \L_2$.  It is enough to prove
the result after extending the base field from $\k$ to $\kbar$, since base
extension of linear maps defined over $\k$ does not change the dimension of 
their kernels or images.  This principle will be used constantly in this
article.

We assume without loss of generality that $\d_1 \leq \d_2$, and begin with
the following exact sequence of vector bundles on $\C$ (using the notation
of \cite{Lazar}):  
\begin{equation}
\label{equation2.3}
0 \to \mathcal{M}_{\L_1} \to \HoLone \tensor \OC \xrightarrow{\alpha}
\L_1 \to 0.
\end{equation}
Here $\alpha$ is the ``evaluation'' map, and $\mathcal{M}_{\L_1} = \ker
\alpha$.  We have used the fact that $\L_1$ is base point free, i.e., that
it is generated by global sections.  Tensoring \eqref{equation2.3} with 
$\L_2$ and taking cohomology reduces our problem to showing that
$H^1(\mathcal{M}_{\L_1} \tensor \L_2) = 0$.  To this end, choose $r = \d_1
- 1 - \g$ points $\P_1, \dots, \P_r$ on $\C$ in general position (this is
where we really need to work over $\kbar$).  Then, by Lemma~1.4.1
and the beginning of the proof of Theorem~1.2.7 of \cite{Lazar}, we obtain
the following exact sequence (where we have temporarily lapsed back into
writing line bundles multiplicatively):\footnote{
	This is where we use the fact that $\d_1 \geq 2\g + 1$, so
	$r \geq \g$, and we can guarantee that $\L_1(-\sum \P_i)$ not only
	has a two-dimensional space of global sections, but is also base
	point free for a generic choice of points $\P_1, \dots \P_r$,
	because it is a generic line bundle of degree $\g+1$.  We then
	obtain \eqref{equation2.4} by combining \eqref{equation2.3} with an
	analogous exact sequence for $\L_1(-\sum \P_i)$, and using the
	snake lemma.  We also note that the analog of \eqref{equation2.3}
	allows us to deduce that $\L_1^{-1}(\sum \P_i) \isomorphic
	\mathcal{M}_{\L_1(-\sum \P_i)}$. 
	This is known as the ``base point free pencil trick''; for
	more details, refer to the exposition in \cite{Lazar} and to
	\cite{ACGH}, pp.~126 and 151, particularly Exercise K. 
	}
%
% note for myself: since there are at least g points, we can ``get anything
% in the Jacobian'' in the sense that LL := L_1(-sum P_i) is any degree g+1
% line bundle.  Now asking LL to have h^0 = 2 and being basepointfree
% amounts to requiring that H^1(LL - Q) = 0 for all Q.  Write K for the
% canonical divisor.  The space of such LLs is that for which the degree
% g-2 bundle K - LL + Q never has sections for any Q, so is never
% equivalent to a sum of g-2 points.  The bad LL are thus those of the 
% form K + (any pt Q) - (any g-2 pts) which are a (g-1) dimensional 
% subspace of the ``full Jacobian'' of degree g+1 divisors LL.  So the
% generic LL is good.  This is the proof of the standard result that the
% generic divisor of degree g+1 is base point free.
%
\begin{equation}
\label{equation2.4}
0 \to \L_1^{-1}(\sum_{i=1}^r \P_i) \to \mathcal{M}_{\L_1} \to 
 \directsum_{i=1}^r \OC (-\P_i) \to 0.
\end{equation}
Tensoring this sequence with $\L_2$ and taking cohomology reduces our
problem to showing, first, that $H^1(\L_2(-\P_i)) = 0$ for all $i$ (as is
clear by using Serre duality and looking at degrees), and, second, that
$H^1(\L_2 \tensor \L_1^{-1}(\sum_{i=1}^r \P_i)) = 0$.  This last statement
is equivalent to requiring $H^0(\mathcal{K} \tensor \L_2^{-1} \tensor \L_1
(- \sum_{i=1}^r \P_i)) = 0$, where we write $\mathcal{K}$ for the canonical
bundle on $\C$.  Since we have assumed that $\d_2 \geq \d_1$, we see that
the line bundle $\L_3 = \mathcal{K} \tensor \L_2^{-1} \tensor \L_1$ has
degree at most $2\g - 2$.  Therefore $\dim H^0(\L_3) \leq \g$, because we
can embed $H^0(\L_3)$ into the $\g$-dimensional space $H^0(\L_3(\D))$,
where $\D$ is any effective divisor of degree $2\g - 1 - \deg \L_3$.  Since 
the $r$ points $\P_1, \dots, \P_r$ impose $r \geq \g$ zeros, we obtain that 
$H^0(\L_3(-\sum_{i=1}^r \P_i)) = 0$ for generic $\P_1, \dots, \P_r$, as
desired.
\end{proof}

Our last result in this section is a converse of sorts to
Lemma~\ref{lemma2.1}, since it describes a kind of ``division'' of sections 
of line bundles.  This relatively easy lemma/algorithm is crucial to our
results, perhaps even more so than the somewhat deeper Lemma~\ref{lemma2.1}.

\begin{lemmalg}[``Division'']
\label{lemmalg2.2}
Let $\C$ be a curve as in Lemma~\ref{lemma2.1}.  Let $\L_1$ and $\L_2$ be
line bundles (written multiplicatively in this lemma/algorithm), and let
$\D_1$ and $\D_2$ be effective divisors on $\C$, where everything is
defined over $\k$.  (We allow $\D_1$ or $\D_2$ to be zero.)  Assume that
$\L_2(-\D_2)$ has no base points (over $\kbar$); for example, it is
sufficient to have $\deg \L_2 - \deg \D_2 \geq 2\g$.  Let
 $\mul: \HoLone \tensor \HoLtwo \to H^0(\L_1 \tensor \L_2)$
be the multiplication map, as before.  Assume that we know the subspaces
$H^0((\L_1  \tensor \L_2)(-\D_1 - \D_2)) \subset H^0(\L_1 \tensor \L_2)$, and
$H^0(\L_2(-\D_2)) \subset \HoLtwo$.  Then we can compute $H^0(\L_1(-\D_1))$ 
to be the space
\begin{equation}
\label{equation2.5}
\{ \s \in \HoLone \mid 	\forall t \in H^0(\L_2(-\D_2)), \quad
	\mul(\s \tensor \t) \in H^0((\L_1 \tensor \L_2) (-\D_1 - \D_2)) \}.
\end{equation}
Note that for computational purposes, it is enough to let $\t$ range only
over a basis for $H^0(\L_2(-\D_2))$ (and not over the whole space) in
equation~\eqref{equation2.5}.
\end{lemmalg}

\begin{proof}
The fact that $H^0(\L_1(-\D_1))$ is a subset of the space defined in
\eqref{equation2.5} is immediate.  As for the 
reverse inclusion, it follows from the fact that $H^0(\L_2(-\D_2))$ has no
base points.  To see why, extend scalars to $\kbar$, and write $\D_1 =
\sum_\P \mP \P$, and $\D_2 = \sum_\P n_\P \P$, where only finitely many of
the $\mP$ and $n_\P$ are nonzero.  For any point $\P$ on $\C$, base point freeness means 
that we can find $\t \in H^0(\L_2(-\D_2))$ which (when viewed as a section
of $\HoLtwo$) vanishes at $\P$ to order exactly $n_\P$.  Thus, if $\s \in
\HoLone$ satisfies the condition in equation~\eqref{equation2.5}, then
$\mul(\s \tensor \t)$ vanishes at $\P$ to order at least $\mP + n_\P$, from
which we conclude that $\s$ vanishes at $\P$ to order at least $\mP$.
Since $\P$ was arbitrary, we obtain that $\s \in H^0(\L_1(-\D_1))$, as
desired.
\end{proof}

\begin{remark}
\label{remark2.3}
Lemma/Algorithm~\ref{lemmalg2.2} works even if the output,
$H^0(\L_1(-\D_1))$, fails to be base point free.  The fact that our section
$\s$ is required to vanish at $\D_1$ may well have severe consequences
for the behavior of $\s$ at other points of $\C$.
\end{remark}

\begin{remark}
\label{remark2.4}
In Remark~\ref{remark3.4.7} below, we shall see that actually calculating
$H^0(\L_1 \tensor \L_2)$ as the image of $\mul$, as in Lemma~\ref{lemma2.1},
can be done in $O(\thedim^4)$ steps, if the dimensions of the vector spaces
involved are of the order of $\thedim$.  (In our applications, $\thedim$ is
a small constant times $\g$.)  We caution the reader to avoid an overly naive
implementation of the calculation that would take time $O(\thedim^5)$.  We 
shall also see in the same remark below that the operation in
Lemma/Algorithm~\ref{lemmalg2.2} can be implemented in time $O(\thedim^4)$.
\end{remark}

\section{Basic algorithms on divisors}
\label{section3}

As in Section~\ref{section2}, let $\C$ be a curve of genus $\g$ over a
perfect field $\k$, and fix a line bundle $\L$ on $\C$ of large degree
$\N$, which will be at least $2\g + 1$, and typically somewhat larger. 
Throughout this section, except briefly in
Theorem/Algorithm~\ref{thmalg3.7}, the letters $\D$, $\D'$, and $\E$ will
refer to effective divisors on $\C$ (that are rational over $\k$), of
degrees $\d$, $\d'$, and $\e$.

We represent an effective divisor $\D$ of degree $\d \leq \N - 2\g$ by the
codimension $\d$ subspace $\WD = \HoLmD$ of the vector space $\V = \HoL$;
alternatively, we can work with its annihilator $\WDstar \subset \Vstar$,
corresponding to the span of the points of $\D$ (counting multiplicity) in
the projective embedding of $\C$ into $\Pr(\V)$.  In later sections, we
give various choices of parameters that allow us to apply the algorithms of
this section to compute in the Jacobian of $\C$.  For ease of following
the exposition in this section, we advise the reader to keep in mind the
application given in Section~\ref{section4}.  In that setting, $g \geq 1$,
and $\L$ will have the form $3\L_0 = \OC(3\D_0)$ for $\D_0$ a divisor of
degree $\d_0 \geq 2\g$; there is some benefit in even assuming that $\d_0
\geq 2\g + 1$, as this allows us to use either Theorem/Algorithm
\ref{thmalg3.4} or \ref{thmalg3.6.7}, whereas if we use $\d_0 = 2\g$, we
can only use the latter Theorem/Algorithm.  The reader should probably use
the value $\d_0 = 2\g + 1$, whence $\N = 3\d_0 = 6\g + 3$.  The values of
$\d$ that will appear in our application are $\d = \d_0 (= 2\g + 1$, say)
and $\d = 2\d_0 (= 4\g + 2)$; an effective divisor $\D$ of degree $\d_0$
will be called {\em small}, and a divisor of degree $2\d_0$ will be called
{\em large}.  A small divisor $\D$ will describe a point on the Jacobian
$\J$ corresponding to the line bundle $\OC(\D - \D_0)$, and a large divisor
$\D'$ will usually arise as an effective divisor equivalent to $3\D_0 - \D$
for some small divisor $\D$.  Nonetheless, we carry through our analysis
for more general values of $\N$ and $\d$, subject to certain inequalities.
This allows for different versions of our algorithms, depending on which
parameters and design choices are most suitable in a given context.

We begin by identifying the subspaces $\WD \intersect \WE$ and $\WD + \WE$,
in terms of divisors related to $\D$ and $\E$.  These subspaces are dual to
$\WDstar + \WEstar$ and $\WDstar \intersect \WEstar$.  From either point of
view (either thinking of $\WD$ as a space of sections which vanish at $\D$,
or of $\WDstar$ as the span of the points of $\D$), it is very plausible
that these new subspaces should correspond to the divisors $\D\union\E$ and
$\D\intersect\E$, at least if all the multiplicities of points are at most
1.  The content of Proposition/Algorithm~\ref{propalg3.2} is that a result
of this kind holds in general, provided we make the following definition.

\begin{definition}
\label{definition3.1}
If $\D = \sum m_\P \P$ and $\E = \sum n_\P \P$ are two (effective)
divisors, then we define
\begin{equation}
\label{equation3.1}
\D \union \E = \sum \max(m_\P, n_\P) \P,
\qquad \qquad
\D \intersect \E = \sum \min(m_\P, n_\P) \P.
\end{equation}
Note that if $\D$ and $\E$ are both defined over $\k$, then so are $\D
\union \E$ and $\D \intersect \E$.
\end{definition}

\begin{propalg}[Union and intersection]
\label{propalg3.2}
Given $\D$ and $\E$ as in Definition~\ref{definition3.1}, then in all cases
$\WD \intersect \WE = \W_{\D \union \E}$ (or, dually, $\WDstar + \WEstar =
\Wstar_{\D \union \E}$).  If furthermore $\deg (\D \union \E) \leq \N -
2\g$, then we can recover $\D \union \E$ from the subspace $\W_{\D \union \E}$,
and we also obtain $\WD + \WE = \W_{\D \intersect \E}$ (dually,
$\WDstar \intersect \WEstar = \Wstar_{\D \intersect \E}$).
\end{propalg}

\begin{proof}
The first statement is trivial, as we are requiring sections to vanish at
$\P$ simultaneously to order at least $m_\P$ and at least $n_\P$.  The
condition on the degree of $\D \union \E$ is to ensure that
$\HoLm{(\D\union\E)}$ is base point free by Lemma~\ref{lemma2.0.5}, so that
we can recover $\D \union \E$.  Moreover, $\codim \W_{\D\union\E} = \deg(
\D \union \E)$ in that case.  As for $\W_{\D \intersect \E}$, we always
have $\WD + \WE \subset \W_{\D\intersect\E}$.  If we know that $\deg(\D
\union \E) \leq \N - 2 \g$, then we can conclude equality by comparing
dimensions, since $\codim (\WD + \WE) = \codim \WD + \codim \WE -
\codim(\WD \intersect \WE)$, and an analogous equality holds for the
degrees of the corresponding divisors.
\end{proof}

\begin{remark}
\label{remark3.1.5}
In the context of Section~\ref{section4},
Proposition/Algorithm~\ref{propalg3.2} allows us to 
find the union and intersection of two small divisors $\D$ and $\E$, since
in that case the degree of $\D \union \E$ cannot exceed $\d + \e = 2 \d_0 = 
\N - \d_0 \leq \N - 2\g$.  In case these divisors are disjoint, as can be
measured by looking at their intersection, this gives a simple way to
compute the sum of these divisors.  We formalize this below, and throw in
an easy algorithm to see if one divisor is contained in the other.  By this
we mean the following definition:
\begin{equation}
\label{equation3.1.5}
\D = \sum m_\P \P \subset \E = \sum n_\P \P \iff \text{ for all }\P, m_\P
\leq n_\P.
\end{equation}
\end{remark}

\begin{propalg}[Disjointness and inclusion]
\label{propalg3.3}
Given divisors $\D$ and $\E$, of degrees $\d$ and $\e$ respectively, then:
\begin{enumerate}
\item If $\d + \e \leq \N - 2\g$, then we can test if $\D \intersect \E =
\emptyset$ by checking if $\WD + \WE = \V$ (alternatively, $\WDstar
\intersect \WEstar = 0$).  In that case, $\D+\E = \D \union \E$ can be
calculated as $\W_{\D + \E} = \WD \intersect \WE$, by
Proposition/Algorithm~\ref{propalg3.2}.  An alternative to checking if $\WD 
+ \WE = \V$ is to simply compute $\WD \intersect \WE$ and to see if the
intersection has codimension $\d + \e$.
\item If $\d, \e \leq \N - 2\g$, then we can test if $\D \subset \E$ (in
the sense of \eqref{equation3.1.5}) by testing whether $\WE \subset \WD$
(dually, if $\WDstar \subset \WEstar$).
\end{enumerate}
\end{propalg}

\begin{remark}
\label{remark3.3.5}
The linear algebra calculations to which we reduce our problem in our
algorithms can be tested easily using standard techniques of Gaussian
elimination.  This involves choosing a basis for $\V$ (hence an isomorphism
$\V \isomorphic \k^{\N + 1 - \g}$), and expressing everything in terms of
the resulting coordinates on $\V$.  A subspace such as $\WD$ can then be
represented in terms of a basis for $\WD$, expressed in terms of the
coordinates on $\V$.  Alternatively, one can work with $\WDstar$ instead of 
$\WD$, since for small divisors, $\WDstar$ has smaller dimension.  The
drawback to this is that using $\WDstar$ appears to be less efficient for our
later algorithms, which involve the multiplication map $\mul$.
We also refer the reader to a side observation in Remark~\ref{remark3.4.7}
for an idea of how to test whether $\D$ and $\E$ are disjoint even if $\d +
\e > N - 2\g$; however, that method requires the slightly stronger
hypothesis that $\d, \e \leq \N - 2\g - 1$.  Alternatively, one can replace 
$\L$ with a line bundle of higher degree, as in the discussion at the end
of Section~\ref{section5}.
\end{remark}

We now describe an algorithm to add two divisors in general, whether or not 
they are disjoint.  We need to assume knowledge of the multiplication map
\begin{equation}
\label{equation3.2}
\mul: \V \tensor \V = \HoL \tensor_\k \HoL \to \HotwoL.
\end{equation}
In the case that interests us, $\N \geq 2\g + 2$.  Then knowledge of $\mul$
is tantamount to knowing the quadratic polynomials that define $\C$ in
the projective embedding given by $\L$, as mentioned just before equation
\eqref{equation2.1}.  The map $\mul$ can be represented explicitly in
terms of bases for $\V$ and for $\HotwoL$.  The technique in the theorem
below is fundamental, and plays a role in all the remaining algorithms of
this article.  We give a geometric interpretation of this technique after 
the proof of Proposition/Algorithm~\ref{propalg3.5} below.  In
Theorem/Algorithm~\ref{thmalg3.6.7} below, we give a second algorithm for
addition of divisors, using some similar ideas.  The other algorithm works
in a slightly more general context than the following algorithm, but the
first algorithm is perhaps simpler to explain.

\begin{thmalg}[Addition of divisors (first method)]
\label{thmalg3.4}
Let $\D$ and $\E$ satisfy $\d, \e \leq \N - 2\g - 1$.  Then we can compute
$\W_{\D+\E}$ in the following steps:
\begin{enumerate}
\item Restricting $\mul$ to the subspace $\WD \tensor \WE = \HoLmD \tensor
\HoLm{\E}$ of $\V \tensor \V$, we obtain from Lemma~\ref{lemma2.1} that
$\mul(\WD\tensor\WE) = \HotwoLm{\D - \E}$.  
%%% note this is H^0(2L-D-E) !!!
The subspace $\HotwoLm{\D - \E}$ of $\HotwoL$ can be computed as the span
of $\mul(\s_i \tensor \s'_j)$, as $\s_i$ and $\s'_j$ range over bases for
$\WD$ and $\WE$ respectively.
\item Using Lemma/Algorithm~\ref{lemmalg2.2}, we compute
\begin{equation}
\label{equation3.3}
\W_{\D+\E} = \{ \s \in \V \mid \forall \t \in \V,
\quad \mul(\s \tensor \t) \in \HotwoLm{\D-\E} \}.
\end{equation}
Recall that one only needs to let $\t$ range over a basis for $\V$.
\end{enumerate}
\end{thmalg}

\begin{proof}
In order to apply Lemma/Algorithm~\ref{lemmalg2.2}, we need to observe that 
$\HoL$ is base point free.  This holds by Lemma~\ref{lemma2.0.5}, since
our assumptions imply that $\N \geq 2\g + 1 +\d \geq 2\g$.
\end{proof}

\begin{remark}
\label{remark3.4.5}
In order for knowledge of $\W_{\D+\E}$ to allow us to uniquely recover
$\D+\E$, we need to know that $\d + \e \leq \N - 2\g$; in the setting of
Section~\ref{section4}, for example, this would require both $\D$ and $\E$
to be small (and $\d_0 \geq 2\g + 1$), even though the above algorithm
correctly calculates $\W_{\D+\E}$ in case one of the divisors is small, and
the other large.  Nonetheless, for use in Theorems/Algorithms
\ref{thmalg4.1} and~\ref{thmalg5.5} and elsewhere, we have presented the
above algorithm in the general case, even if $\W_{\D+\E}$ does not uniquely
determine $\D+\E$.  We note as an aside that this last case happens
precisely when $\HoLm{\D-\E}$ has base points, and in that case the
inclusion
\begin{equation}
\label{equation3.4}
\mul(\HoLm{\D-\E} \tensor \HoL) \subset \HotwoLm{\D-\E}.
\end{equation}
will be strict, since $\HotwoLm{\D-\E}$ will never have base points, by
Lemma~\ref{lemma2.0.5}.  (This of course does not affect the validity of our 
proof; in the easy half of the proof of Lemma/Algorithm~\ref{lemmalg2.2}, 
we only need the inclusion in equation~\eqref{equation3.4} to obtain that
$\W_{\D+\E} \subset \{ \s $ as above$\}$.)
\end{remark}

\begin{remark}
\label{remark3.4.7}
We discuss some practical issues in implementing the above algorithm,
beyond the fact that Theorem/Algorithm~\ref{thmalg3.4} settles the
theoretical issue of computing $\W_{\D+\E}$.  We also discuss the running
time of the algorithm.  Let the dimensions of the vector spaces involved in
the calculation be of the order of $\thedim$ (this is approximately $\dim V
= \N + 1 - \g$, which is $O(g)$ for all our applications).  Then an overly
naive implementation of the algorithm as stated above takes time
$O(\thedim^5)$, but a little care brings this down to $O(\thedim^4)$.  We
can use some randomness to make the algorithm slightly faster, but still
taking $O(\thedim^4)$ steps.  (We remind the reader that we count an
operation in the field $\k$ as a single step, taking one unit of time.)  We
note that our second algorithm below, Theorem/Algorithm~\ref{thmalg3.6.7},
also takes time $O(\thedim^4)$.

We first discuss why the naive implementation of
Theorem/Algorithm~\ref{thmalg3.4} runs in time $O(\thedim^5)$, and how to
reduce this time to $O(\thedim^4)$.  We assume that $\mul$ is implemented
in terms of fixing a basis $\{\t_i\}$ for $\V$, and storing all the
products $v_{ij} = \mul(\t_i \tensor \t_j) \in \HotwoL$.  Moreover, a
subspace $\WD$ is represented by choosing a basis, and writing all the
basis elements in terms of the $\t_i$.  Then it takes $O(\thedim^3)$
operations to calculate each individual $\mul(\s_i \tensor \s'_j)$, as
$\s_i$ and $\s'_j$ range over bases for $\WD$ and $\WE$.  Since there are
$O(\thedim^2)$ such pairs $(\s_i, \s'_j)$, it seems that we will need
$O(\thedim^5)$ steps merely to write down a spanning set for
$\HotwoLm{\D-\E}$, viewed as a subspace of the $O(\thedim)$-dimensional
space $\HotwoL$.  However, we can proceed more quickly, at the expense of
using more memory, as follows: first calculate (and store) each product
$w_{ij} = \mul(\s_i \tensor \t_j)$, as $\s_i$ ranges over a basis for $\WD$,
and $\t_j$ ranges over the standard
basis for $\V$.  This takes only $O(\thedim^4)$
steps, since each individual multiplication $\mul(\s_i \tensor \t_j)$ now
needs only $O(\thedim^2)$ operations, as it involves only the $v_{i'j}$
for a fixed $j$ and varying $i'$.  Now if $\E = 0$, we are already done;
more generally, we use the $w_{ij}$ to assemble all the products $\mul(\s_i
\tensor \s'_j)$ in a further $O(\thedim^4)$ steps.  Since these products
will virtually always be linearly dependent, we should then use Gaussian
elimination to reduce these products to a basis for $\HotwoLm{\D - \E}$.
This takes $O(\thedim^4)$ operations, since we have $O(\thedim^2)$ vectors
to reduce inside the $O(\thedim)$-dimensional space $\HotwoL$.\footnote{It
would be interesting to try to use the fact that for fixed $i$ and varying
$j$ (or vice-versa), the products $\mul(\s_i \tensor \s'_j)$ are known to
be linearly independent before we start Gaussian elimination.  Still, there
are advantages to obtaining a basis for $\HotwoLm{\D-\E}$ in echelon form,
once we invoke Lemma/Algorithm~\ref{lemmalg2.2}.}

All of this takes some time.  So in practice (at least when $\d + \e
\leq \N - 2\g$), one should rely on the fact that the divisors $\D$ and
$\E$ are likely to be disjoint, and simply compute $\WD \intersect \WE$, as
mentioned in Proposition/Algorithm~\ref{propalg3.3}.  This only takes
$O(\thedim^3)$ operations.  If the intersection has the right codimension
$\d + \e$, which should happen rather often, then we immediately obtain
$\W_{\D + \E}$ as this intersection.  We note that in the undesirable case
$\d + \e > \N - 2\g$, it does not seem to be worth the computational effort
to check if $\D$ and $\E$ are disjoint; one might as well use
Theorem/Algorithm~\ref{thmalg3.4} directly.  (Here is one idea of how one
can adapt our algorithms to check that $\D$ and $\E$ are disjoint, even if
$\d+\e > \N - 2\g$: first calculate $\HotwoLm{\D} = \mul(\WD
\tensor \V)$, and similarly $\HotwoLm{\E}$; this works, and takes time
$O(\thedim^4)$ as before.  Then find the intersection $\HotwoLm{\D}
\intersect \HotwoLm{\E}$ inside the larger ambient space $\HotwoL$ instead
of $\HoL$.  This amounts to using $2\N$ instead of $\N$, which allows us to
proceed as in Proposition/Algorithm~\ref{propalg3.3}, since $\d + \e$ is
now less than $2\N - 2\g$.)

Coming back to the issue of calculating $\mul(\WD \tensor \WE)$, we note
that we can also try to avoid calculating all $O(\thedim^2)$ products
$\mul(\s_i \tensor \s'_j)$.  Instead, we randomly choose $O(\thedim)$ pairs 
of the form $(\s, \s')$, where $\s \in \WD$ and $\s' \in \WE$ are random
(nonzero) elements.  In that case, the resulting products $\mul(\s \tensor
\s')$ are likely to already span $\HotwoLm{\D - \E}$, especially 
if we take, say, twice as many pairs $(\s,\s')$ as the dimension of
$\HotwoLm{\D - \E}$.  With high probability, this reduces the time we spend 
on finding $\mul(\WD \tensor \WE)$, but it still involves $O(\thedim^4)$
operations, since one must calculate the $\mul(\s \tensor \s')$.  On the
other hand, this reduces the time spent on Gaussian elimination to
$O(\thedim^3)$, but without an assurance of success as we have not
necessarily used a full spanning set for $\HotwoLm{\D - \E}$.

As for the ``division'' calculation in step~2 of
Theorem/Algorithm~\ref{thmalg3.4}, it can be done in $O(\thedim^4)$
operations as follows: first, find a matrix $M$, representing a linear
transformation from $\HotwoL$ to $\k^{\d + \e}$, whose kernel is precisely
$\HotwoLm{\D-\E}$.  This amounts to finding a basis of the annihilator
$\HotwoLm{\D-\E}^* \subset \HotwoL^*$, and can be done in time at most
$O(\thedim^3)$ (even $O(\thedim^2)$, if the basis for $\HotwoLm{\D-\E}$ is
already in echelon form).  Recall that we have on hand all the products
$v_{ij} = \mul(\t_i \tensor \t_j) \in \HotwoL$, where the $\t_i$ are a
basis for $\V$.  We now need to solve for $\s = \sum_i x_i \t_i$ such that
for all $j$, $\sum_i x_i M(v_{ij}) = (0, \dots, 0) \in \k^{\d + \e}$.
Computing the various $M(v_{ij})$ takes $O(\thedim^4)$ steps, and the
resulting Gaussian elimination to solve for the tuples $(\dots, x_i,
\dots)$ representing $\s$ takes $O(\thedim^4)$ operations, being a system
of $O(\thedim^2)$ equations in $O(\thedim)$ variables.  We point out that
we will later still be able to perform similar calculations (inspired by
Lemma/Algorithm~\ref{lemmalg2.2}) in $O(\thedim^4)$ operations, even if we
no longer have $\t$ ranging over all of $\V$.  This is necessary for 
equation~\eqref{equation3.6} in Proposition/Algorithm~\ref{propalg3.5}, as
well as Theorems/Algorithms \ref{thmalg3.6}, \ref{thmalg3.6.7},
and~\ref{thmalg3.7}.  The idea is that we need to compute $v'_{ij} =
\mul(\t_i \tensor \t'_j)$ for $\t_i$ as before, but with 
$\t'_j$ ranging over a basis for some suitable subspace of $\V$.  Each
$v'_{ij}$ can be computed from the existing collection of $v_{ij}$ in
$O(\thedim^2)$ operations, so it takes time $O(\thedim^4)$ to assemble the
$v'_{ij}$, and the calculation proceeds as before.
\end{remark}

Our next algorithm is a straightforward illustration of the ideas in
Theorem/Algorithm~\ref{thmalg3.4}.  Given two divisors $\D = \sum m_\P \P$
and $\E = \sum n_\P \P$, we may wish to calculate the divisor
\begin{equation}
\label{equation3.5}
\E \, \backslash \, \D = \sum \max(0, n_\P - m_\P) \P.
\end{equation}
In case the multiplicities are at most 1, then this is the usual operation
on sets of points.  If $\D \subset \E$, this is just subtraction of
divisors.

\begin{propalg}[Set subtraction]
\label{propalg3.5}
Given $\D$ and $\E$, with $\d \leq \N - 2\g$ and $\e \leq \N - 2\g - 1$, we
can calculate $\W_{\E \backslash \D}$ as follows: 
First, calculate $\HotwoLm{\E} = \mul(\W_{\E} \tensor \V)$.  Then compute
\begin{equation}
\label{equation3.6}
\W_{\E \backslash \D} = \{ \s \in \V \mid \forall \t \in \WD, \quad
	\mul(\s \tensor \t) \in \HotwoLm{\E} \}.
\end{equation}
\end{propalg}

\begin{proof}
Lemma~\ref{lemma2.1} gives us the surjectivity of $\mul: \HoLm{\E} \tensor
\HoL \to \HotwoLm{\E}$.  As for \eqref{equation3.6}, it follows from a
modification of the argument in Lemma/Algorithm~\ref{lemmalg2.2}, since
$\HoLmD$ has no base  points.  Note again that the computation on
\eqref{equation3.6} can be done effectively by letting $\t$ range only over
a basis for $\WD$.
\end{proof}

We now come to the most important basic algorithm in this section.  It is
close to Theorem/Algorithm~\ref{thmalg3.4} and
Proposition/Algorithm~\ref{propalg3.5}, but we shall first pause to
describe what the algorithm is doing geometrically.  Our discussion will
also shed light on Theorem/Algorithm~\ref{thmalg3.4}.  The geometric
picture is best understood in terms of the projective embedding of $\C$
into $\Pr(\V)$ given by $\L$, and in terms of ideals of the graded
projective coordinate ring of $\Pr(\V)$.  To this end, let $\A =
\Sym^* \V$ be this projective coordinate ring; $\A$ is the symmetric
algebra on $\V$, and an element of $\V$ is a linear (degree one) polynomial
whose vanishing set defines a hyperplane in $\Pr(\V)$.  The curve $\C$ is
defined by a homogeneous ideal $\I_\C$, which is generated in our
application by elements of degree two, and these elements can more or less
be identified with the kernel of $\mul : \V\tensor\V \to \HotwoL$.  (The
quotient $\A / \I_\C$ is the projective coordinate ring of $\C$, and $\mul$
is just multiplication in this quotient.)  Given a divisor $\D$ and its
associated subspaces $\WD$ and $\WDstar$, we shall abuse notation and also
write $\WDstar$ for the corresponding plane in $\Pr(\V)$, as described just
after Lemma~\ref{lemma2.0.5}.  In that case, we write $\I_\W$ for the ideal
defining this plane; $\I_\W$ is then generated by $\WD$, which consists of
the linear functions vanishing on $\WDstar$.  As we have seen before, $\D =
\WDstar \intersect \C$.  In terms of ideals, this means that $\D$ is
defined by the ideal $\I_\D = \I_\W + \I_\C = \WD \A + \I_\C$.

Now our next algorithm implements the following operation: given a section
$\f \in \WD$, we know that $\f$ (viewed as a section of $\L$) vanishes at a 
divisor of the form $\E = \D + \D'$ belonging to the equivalence class of
$\L$.  We then wish to find $\W_{\D'}$, essentially by the same method as in
Proposition/Algorithm~\ref{propalg3.5}.  Geometrically, the zero set of
$\f$ is a hyperplane $\Hyp$ passing through $\D$, and its ideal $\I_\Hyp$
is generated by $\f$.  Then the divisor $\E$ is the intersection $\Hyp
\intersect \C$, and is defined by $\I_\E = \I_\Hyp + \I_\C = \f \A +
\I_\C$.  As a zero-dimensional scheme, we then obtain that $\D' = \E \,
\backslash \, \D$ should be defined by 
\begin{equation}
\label{equation3.7}
\I_{\D'} = (\I_\E : \I_\D) = \{ a \in \A \mid a \I_\D \subset \I_\E \} =
\{a \in \A \mid a \WD \subset \f \A + \I_\C \}.
\end{equation}
Namely, functions should vanish on $\D'$ if and only if when we multiply
them by a function vanishing on $\D$, we obtain a function vanishing on $\E
= \D \intersect \D' = \Hyp \intersect \C$.  This fact is not so immediate
if $\D$ and $\E$ have multiplicities greater than~1; it is essential, for
example, that our divisors lie on a fixed curve, so that every extra order
of vanishing at a point $\P$ imposes exactly one more condition on
homogeneous polynomial functions, provided that the degree of $\E$ is not
too large.  At any rate, we are not interested in the ideal defining $\D'$
itself, but rather in the plane $\Wstar_{\D'}$ spanned by $\D'$.  Dually,
we only wish to determine $\W_{\D'}$, which consists of the linear
functions vanishing on $\D'$.  Thus we only seem to need to find the degree
one elements of $\I_{\D'}$.  Using~\eqref{equation3.7}, this leads us to
the calculation that we present below in Theorem/Algorithm~\ref{thmalg3.6}.
(Multiplication of linear functions and taking the quotient by the
quadratic elements of $\I_\C$ is exactly the map $\mul$.)  Again, the fact
that we can content ourselves with the degree one elements of $\I_{\D'}$ is
not so trivial: in general, the vanishing set of a homogeneous ideal on
projective space depends only on the saturation of that ideal (exercise
II.5.10 of \cite{Harts}; note that step~2 in
Theorem/Algorithm~\ref{thmalg3.4} is effectively computing a saturation).
Now knowing the degree one elements of a saturation can in principle
involve elements of arbitrarily high degree in the ideal.  To do this
algorithmically in general would probably involve Gr\"obner basis
techniques, which will typically be slower than our methods --- the
Gr\"obner basis calculations can take exponential time in the worst case
(although they are efficient for ``most'' calculations), and even their
average case behavior is likely to be slower than our linear algebra
calculations.  The point of our techniques is that we can use them to
justify the geometric reasoning in the above heuristic argument.  We
therefore (finally) present the following result.

\begin{thmalg}[Flipping algorithm]
\label{thmalg3.6}
Let $\D$ be an effective divisor of degree $\d \leq \N - 2\g$.  Take a
nonzero element $\f \in \WD$, and write the divisor of $\f$ (viewed as a
section of $\L$) as
\begin{equation}
\label{equation3.8}
(\f) = \E = \D + \D'.
\end{equation}
Thus $\d'= \deg \D' = \N - \d$, and $\D'$ belongs to the divisor class of
$\L - \D$.  Then we can compute $\W_{\D'}$ in the following steps:
\begin{enumerate}
\item
Calculate $\HotwoLm{\D-\D'} = \mul(\f \tensor \HoL)$.
\item
Compute $\W_{\D'} = \{ \s \in \HoL \mid$ for all $\t$ in (a basis for)
$\WD$, $\mul(\s \tensor \t) \in \HotwoLm{\D - \D'} \}$. 
\end{enumerate}
Note that if furthermore $\d \geq 2\g$, then knowledge of $\W_{\D'}$
determines $\D'$ uniquely, by Lemma~\ref{lemma2.0.5}.
\end{thmalg}

\begin{proof}
The identity in the first step (which amounts to computing the image of the
degree two elements of $\f \A$ in the quotient $\A / \I_\C$) does not need
Lemma~\ref{lemma2.1}, but turns out to be much simpler.  In fact,
multiplication by $\f$ is a bijection between $\HoL$ and
$\HotwoLm{\D-\D'}$, the inverse being division by $\f$ (which does not
introduce any poles!).  Then the second step above (which is the degree one 
part of the left hand side of \eqref{equation3.7}, whose computation
involves some degree two elements on the right hand side) follows from 
Lemma/Algorithm~\ref{lemmalg2.2}, since $\HoLmD$ has no base points.
\end{proof}

\begin{remark}
\label{remark3.6.5}
In the setting of our application in Section~\ref{section4}, we have $\L =
\OC(3\D_0)$, $\N = 3\d_0$, and $\d_0 \geq 2\g$.  The above algorithm
then allows us to go from a divisor $\D$ (which may be small or large,
corresponding to $\d = \d_0$ or $2\d_0$ respectively) to a complementary
divisor $\D'$ (which is then respectively large or small), such that $\D +
\D'$ is linearly equivalent to $3\D_0$.  This is reminiscent of the theory
behind the Weierstrass embedding of an elliptic curve via a cubic equation,
except that we are contenting ourselves with quadratic equations and must
therefore use a more positive line bundle for our embedding of $\C$.  In
Section~\ref{section5}, we go through analogs of our construction using
smaller divisors (e.g., $\d_0 = \g + 1$, and $\N = \deg \L = 3\g + 3$),
working instead with ``higher degree equations,'' roughly in the sense of
multiplications $\V \tensor \V \tensor \dots \tensor \V \to H^0(n\L)$.
\end{remark}

\begin{remark}
\label{remark3.6.6}
Theorem/Algorithm~\ref{thmalg3.6} also runs in time $O(\thedim^4)$, in the
notation of Remark~\ref{remark3.4.7}.  This is because the
``multiplication'' involved in computing $\mul(\f \tensor \V)$ takes
$O(\thedim^3)$ steps, including reducing the basis $\{\mul(\f \tensor
\t_i)\}$ for $\HotwoLm{\D - \D'}$ to echelon form.  (We already know that
the vectors $\{\mul(\f \tensor \t_i)\}$ are linearly independent.)  It then
takes another $O(\thedim^4)$ steps to implement the second ``division''
step of Theorem/Algorithm~\ref{thmalg3.6}, by the same discussion as in
Remark~\ref{remark3.4.7}.
\end{remark}

Now that we have discussed the flipping algorithm, we use a slight variant
of it to give a second algorithm to compute the sum of two divisors.  It is
not clear whether this is more efficient than
Theorem/Algorithm~\ref{thmalg3.4}, as the new algorithm still requires
time $O(\thedim^4)$, although it involves slightly different linear algebra.
On
the other hand, it has slightly weaker hypotheses, and can work in the
setting of Section~\ref{section4} with $\d_0 = 2\g$.

\begin{thmalg} [Addition of divisors (second method)]
\label{thmalg3.6.7}
Let $\D$ and $\E$ be two divisors with $2\g \leq \d \leq \N - 2\g$.  (We 
do not make any assumptions on $\e$.)  Then we can compute $\W_{\D+\E} =
\HoLm{\D - \E}$ as follows:
\begin{enumerate}
\item
Choose a nonzero $\f \in \WD$, whose divisor (viewing $\f$ as a section of
$\L$) is $(\f) = \D + \D'$; proceed as in Theorem/Algorithm~\ref{thmalg3.6}
to compute $\W_{\D'}$.
\item
Also compute $\HotwoLm{\D-\D'-\E} = \mul(\f \tensor \WE)$.
\item
Now compute $\W_{\D+\E} = \{ \s \in \V \mid \forall \t \in \W_{\D'}, \quad
\mul(\s \tensor \t) \in \HotwoLm{\D - \D' - \E} \}$.
\end{enumerate}
\end{thmalg}
\begin{proof}
This is an almost verbatim adaptation of the argument in
Theorem/Algorithm~\ref{thmalg3.6}.  The hypothesis on $\d$ ensures that
$\WD$ and $\W_{\D'}$ have no base points.  Note that, as in
Theorem/Algorithm~\ref{thmalg3.4}, we cannot hope to recover the actual
divisor $\D + \E$ unless $\d + \e \leq \N - 2\g$.
\end{proof}

We conclude this section with an algorithm which tells whether a given
codimension $\d$ subspace $\W \subset \V$ actually is of the form $\WD$ for 
a divisor $\D$.  In other words, given a point on the Grassmannian
$\Gd(\V)$, we wish to identify whether or not it lies in the image of
$\SymdC$.  Since we have not explicitly written down equations describing
this subvariety of the Grassmannian, we have to do something cleverer than
checking whether certain polynomial functions vanish.  The answer turns out 
to be relatively straightforward, using Theorem/Algorithm~\ref{thmalg3.6}.
Once again, the algorithm takes time $O(\thedim^4)$.

\begin{thmalg}[Membership test]
\label{thmalg3.7}
Let $\W \subset \V = \HoL$ have codimension $\d$, where $2\g \leq \d \leq
\N - 2\g$.  (In particular, $\N \geq 4\g$.)  Choose any nonzero $\f \in
\W$, and compute $\W'$ in the same way as in
Theorem/Algorithm~\ref{thmalg3.6}; i.e., 
\begin{equation}
\label{equation3.9}
\W' = \{ \s \in \V \mid \forall \t \in \W, \quad \mul(\s \tensor \t) \in
\mul(\f \tensor \V) \}.
\end{equation}
Then $\W$ is of the form $\WD$ for some (uniquely determined) $\D$ if and
only if the codimension of $\W'$ in $\V$ is $\N - \d$, or equivalently if
$\dim \W' = \d + 1 - \g$.
\end{thmalg}

\begin{proof}
If $\W$ is of the form $\WD$, then the resulting $\W'$ is equal to the
space $\W_{\D'}$ in Theorem/Algorithm~\ref{thmalg3.6}, where the degree of
$\D'$ is $\d' = \N - \d$.  Then $\L - \D'$ has degree $\d$, and the result
follows.  Conversely, given any $\W$, define a divisor $\D = \sum \mP \P$
by the property that for each $\P$, $\mP$ is the minimum order to which any 
nonzero section $\s \in \W$ vanishes at $\P$ (viewing $\s$ as a section of
$\L$).  In other words, $\D$ is the largest possible divisor for which $\W 
\subset \WD$, and $\D$ is the only possible candidate for a divisor which
might give rise to $\W$.  (One easily checks that $\D$ is defined over $\k$,
even if the individual points are only defined over $\kbar$, but we do not
actually need this fact for our proof, and give our proof by working
entirely over $\kbar$.  As usual, there is no problem in extending
scalars from $\k$ to $\kbar$ for the proof.)
Note that we do not yet know whether $\deg D = \d$, in contrast to our
previous notational convention.  We do know, however,
that $\deg D \leq \d$, since otherwise, by Lemma~\ref{lemma2.0.5}, the
codimension of $\WD$ would be at least $\d+1$, contradicting the inclusion
$\W \subset \WD$.  We see in fact that $\deg \D = \codim \WD \leq \d =
\codim \W$, and that $\W = \WD$ if and only if $\deg \D = \d$.

Now our choice of $\f$ gives us in particular a nonzero element of $\WD$,
so the divisor of $\f$ is of the form $\D + \D'$, where 
\begin{equation}
\label{equation3.10}
\d' = \deg \D' = \N - \deg \D \geq \N - \d.
\end{equation}
The same reasoning as in Lemma/Algorithm~\ref{lemmalg2.2} then shows us
that $\W' = \W_{\D'}$ (the main point is that $\mul(\f \tensor \V) =
\HotwoLm{\D - \D'}$, just as in Theorem/Algorithm~\ref{thmalg3.6}; on the
other hand, for each $\P$, we can already find a section $\t \in \W$ 
that vanishes at $\P$ precisely to the minimum possible order $\mP$).  Now
suppose that $\deg \D'$ were different from $\N - \d$.  In that case,
equation~\eqref{equation3.10} would imply that $\d' \geq \N - \d + 1$.  Let
$\E$ be the sum of any $\N - \d + 1$ points of $\D'$; then the codimension
of $\W_{\D'}$ in $\V$ would be greater or equal to the codimension of
$\W_\E$, which we know to be $\N - \d + 1$ (we are using the fact that $\N
- \d + 1 \leq \N - 2\g + 1$, so as to apply Lemma~\ref{lemma2.0.5}).  This
would contradict our assumption on the codimension of $\W_{\D'}$; hence we
conclude that $\deg \D' = \N - \d$ after all, and hence by
\eqref{equation3.10} we conclude that $\deg D = \d$ and that $\W = \WD$, as 
desired.
\end{proof}

\section{The large model: Weierstrass-style algorithms on divisor classes
and on the Jacobian}
\label{section4}

We now describe a way to work with divisor classes on $\C$, which is
inspired by the usual geometric definition of the group law on an elliptic
curve represented by the plane cubic equation $y^2 = x^3 + ax + b$.
Actually, in many cases, we give two algorithms for the basic operations,
the first of which is straightforward and elegant, and the second of which
is somewhat more efficient, and is given to illustrate various techniques.
In Section~\ref{section5}, we give slightly faster algorithms, using a
somewhat more compact representation of divisor classes.  In all cases, the
various algorithms require $O(\thedim^4)$ operations in the field $\k$,
where $\thedim = \dim \V$ is approximately $5\g$.  The main improvement in
speed in Section~\ref{section5} comes from our bringing $\thedim$ down to a
smaller multiple of $\g$ there.  Thus the algorithms in
Section~\ref{section5} supersede the ones here, at the expense of elegance.
The basic structure of all the algorithms is roughly the same in all
settings.  We refer to the setup in this section as the {\em large model}.

We write $\D \linequiv \E$ to denote that the divisors $\D$ and $\E$ are
linearly equivalent, and we denote the equivalence class of $\D$ by $[\D]$.
We keep the notation of Section~\ref{section3}, and specialize the
situation to our application as follows.  Take an effective divisor $\D_0$
of degree $\d_0 \geq 2\g$ (or $2\g + 1$, if we wish to repeatedly use
Theorem/Algorithm \ref{thmalg3.4} instead of~\ref{thmalg3.6.7}).  Fix
the ambient line bundle $\L = \OC(3\D_0)$; thus $\N = 3 \d_0
\geq 6\g$.  Not surprisingly, we also assume that $\g \geq 1$ (but
everything works with $\g=0$ and $\d_0 \geq 1$, if one is desperate or
perverse enough to use these algorithms in that case).  We note that
it may be tricky to find $\D_0$ of a specific degree, unless the curve $\C$
contains one or more explicitly known rational points.  In the best case,
$\d_0 = 2\g$, and the dimension of the space $\V$ is $\thedim = \N + 1 - \g
= 5\g + 1$.

An element of the Jacobian of $\C$ is a linear equivalence class of degree
zero divisors; these can all be obtained as classes of the form $[\D -
\D_0]$, where $\D$ is an effective divisor of degree $\d_0$, called a {\it
small} divisor.  Given such a divisor, it corresponds to the element $\xD$
of the Jacobian, given by
\begin{equation}
\label{equation4.0.2}
\xD =  [\D - \D_0], \quad \text{which is represented by the subspace } \WD
\subset \V.
\end{equation}
(Thus, if $\d_0 = 2\g$, then $\WD$ is a
$(3\g+1)$-dimensional subspace of $\V$.)  We shall also occasionally make
use of {\it large} divisors, namely, divisors of degree $2\d_0$, for
intermediate results in computations.  We have already addressed in
Theorem/Algorithm~\ref{thmalg3.7} the question of how to test whether a
subspace $\W \subset \V$ actually corresponds to an element of the
Jacobian.  We now turn to the question of deciding when two such subspaces,
of the form $\WD$ and $\WE$, represent the same element of the Jacobian.
This amounts to testing whether the divisors $\D$ and $\E$ are linearly
equivalent.

\begin{thmalg}[Equality of divisor classes]
\label{thmalg4.1}
Let $\D$ and $\E$ be divisors of the same degree $\d \leq \N - 2\g$;
typically, they will be either both small, or both large divisors,
represented by $\WD$ and $\WE$.  Choose any nonzero section $\f \in \WD$,
whose divisor (viewing $\f$ as a section of $\L$) is $(\f) = \D + \D'
\linequiv \L \linequiv 3\D_0$.  Then compute $\W_{\D' + \E}$ in {\em one}
of the following two ways:
\begin{itemize}
\item
Use
%the method of
Theorem/Algorithm~\ref{thmalg3.6} to compute the space $\W_{D'}$, and then
compute $\W_{\D' + \E}$ by Theorem/Algorithm \ref{thmalg3.4} or
\ref{thmalg3.6.7}.  This is less efficient but easier to explain than the
method proposed below.  However, this first method further requires
the assumption that $2 \g + 1 \leq \d \leq \N - 2\g - 1$, in case we use
Theorem/Algorithm~\ref{thmalg3.4}.  If we use
Theorem/Algorithm~\ref{thmalg3.6.7} instead, then we can weaken the
assumption to $2\g \leq \d \leq \N - 2\g$.
\item
Alternatively, adapt the method of Theorem/Algorithm~\ref{thmalg3.6.7},
without explicitly calculating $\W_{\D'}$.  Namely, calculate
$\HotwoLm{\D-\D'-\E} = \mul(\f \tensor \WE)$, and then obtain
$\W_{\D'+\E} = \{ \s \in \V \mid \forall \t \in \W_{\D},
\quad \mul(\s \tensor \t) \in \HotwoLm{\D - \D' - \E} \}$.
\end{itemize}
Then, in either case, $\E \linequiv \D$ (i.e., $[\E]=[\D]$, meaning $\xE =
\xD$ on the Jacobian), if and only if the space $\W_{\D' + \E}$ is not
zero, in which case it must be one-dimensional.
\end{thmalg}
\begin{proof}
The second method to calculate $\W_{\D' + \E}$ works by the same reasoning
as in the proof of Theorem/Algorithm~\ref{thmalg3.6.7}.  Now if there
exists a nonzero $\s \in \W_{\D' + \E}$, then its divisor of zeros $(\s)$
(viewing $\s$ as a section of $\L$) is an effective divisor of degree $\N$
which contains $\D' + \E$.  By comparing degrees, we see that $\D' + \E =
(\s)$.  However, $(\s)$ and $(\f)$ are linearly equivalent (both belong to
the class of $\L$), from which we obtain $\E \linequiv \D$.  Conversely,
if we have $\E \linequiv \D$, then $\D' + \E \linequiv \L$, and it
follows that $\W_{\D'+\E} = \HoOC{\L - \D' - \E} \isomorphic H^0(\OC) =
\k$.
\end{proof}

\begin{remark}
\label{remark4.2}
In the first method above, instead of using Theorem/Algorithm
\ref{thmalg3.4} or~\ref{thmalg3.6.7} to compute $\W_{\D' + \E}$, we can
choose the section $\f$ randomly, in which there is a high probability that
the divisors $\D'$ and $\E$ are disjoint.  In that case, we simply check if
there are any nonzero elements in $\W_{\D'} \intersect \WE = \W_{\D' +
\E}$.  As mentioned in Remark~\ref{remark3.4.7}, it may not be worth our
while to test whether $\D'$ and $\E$ are disjoint; however, if the
intersection is zero, then we can at least immediately conclude that $\xE
\neq \xD$, which saves a good amount of time.  Another point to consider is
that if we wish to test the equality of many different classes $[\E]$
against $[\D]$ by the first method, we should probably compute $\W_{\D'}$
just once, instead of choosing a random $\f$ each time.  We note
incidentally that the space $\W_{\D' + \E}$ that we compute has codimension
less than the degree of $\D' + \E$, a fact that follows immediately from
the fact that $\W_{\D'+\E}$ is either zero or one-dimensional.  We refer
the reader to Remark~\ref{remark3.4.5} for some further discussion.
\end{remark}

We now describe the main building block in our algorithms to effectively
compute in the Jacobian of $\C$.  This algorithm is analogous to the main
geometric construction on the Weierstrass model of an elliptic curve in the 
plane.  Namely, given two points $P_1, P_2$ on the elliptic curve, we
draw the line through them and find its third point of intersection $Q$
with the curve.  This means that $[Q] = - ([P_1] + [P_2])$ in the group
law on the elliptic curve.  We generalize this to the case of divisors of
degree $\d_0$ on our curve $\C$.

\begin{propalg}[``Addflip'']
\label{propalg4.2.5}
Consider two subspaces $\W_{\D_1}, \W_{\D_2}$ associated to two small
divisors $\D_1, \D_2$ on $\C$, corresponding to the classes $\x_{\D_1}$,
$\x_{\D_2}$ on the Jacobian.  The following algorithm then computes the
subspace $\WE$ associated to a divisor $\E$ such that $\xE = -(\x_{\D_1} +
\x_{\D_2})$ (i.e., $\x_{\D_1} + \x_{\D_2} + \xE = 0$):
\begin{enumerate}
\item
First compute $\W_{\D_1 + \D_2}$, either by Theorem/Algorithm
\ref{thmalg3.4} or~\ref{thmalg3.6.7}, or if possible by part~1 of
Proposition/Algorithm~\ref{propalg3.3}.
\item
Now apply Theorem/Algorithm~\ref{thmalg3.6} to obtain the space $\WE$
corresponding to a ``flip'' $\E$ of $\D_1 + \D_2$, where $\D_1 + \D_2 +
\E \linequiv 3\D_0$.
\end{enumerate}
\end{propalg}
\begin{proof}
Immediate.  Note that step~1 is analogous to drawing the line through two
points on a Weierstrass model elliptic curve in the plane, and step~2 is
analogous to finding the third point on the intersection of the elliptic
curve with the line.  Also note that in step~1, our use of
Theorem/Algorithm \ref{thmalg3.4} or~\ref{thmalg3.6.7} means that we
do not need to give special treatment to the case where $\D_1 = \D_2$.
In the Weierstrass model, this would mean that our line through a double
point is automatically calculated as the appropriate tangent line.  It is
however simpler in the Weierstrass model to find the line through two
distinct points.  The analogous statement in our situation is: if 
we know beforehand that $\D_1$ and $\D_2$ are disjoint, 
then we can avoid using the more complicated theorems/algorithms, and
instead simply calculate $\W_{\D_1 + \D_2} = \W_{\D_1} \intersect
\W_{\D_2}$ as in Proposition/Algorithm~\ref{propalg3.3}.  This
will save a lot of time.  In practice, the divisors $\D_1$ and $\D_2$ are
likely to be disjoint, so one should first calculate the intersection, and
only use the more complicated algorithms for adding divisors in case
the codimension of $\W_{\D_1} \intersect \W_{\D_2}$ is not $2\d_0$.
\end{proof}

Now that we have implemented the ``addflip'' operation $f(\x_1,\x_2) =
-(\x_1 + \x_2)$ on the Jacobian of $\C$, we can describe algorithms for all
the group operations.  We begin with negation.  Given a small divisor $\D$,
we seek $\E$ such that $\xE = - \xD$, or in other words $\D + \E \linequiv
2\D_0$.  This can be done in a straightforward way, but we give a slight
modification which uses vector spaces of smaller dimension.  This will make
the negation algorithm more efficient, a point especially worth noting,
since Theorem/Algorithm~\ref{thmalg4.4} (which adds two points on the
Jacobian), involves a negation.  Similar ideas also come into play for
subtraction in the Jacobian, as discussed in
Proposition/Algorithm~\ref{propalg4.5}.

\begin{thmalg}[Negation]
\label{thmalg4.3}
Let $\WD$ be given, where $\D$ is a small divisor.  Then we can compute
$\WE$, where $\xE = - \xD$, in {\em either one} of the following two ways:
\begin{itemize}
\item
Apply Proposition/Algorithm~\ref{propalg4.2.5} to $\D_1 = \D$ and $\D_2 =
\D_0$.  We of course keep the space $\W_{\D_0}$ at hand; it corresponds to
the zero element on the Jacobian.
\item
Alternatively, if one is willing 
to hold in memory more extensive information on the multiplication map
$\mul$, one can find $\WE$ using vector spaces of smaller dimension.
Namely, we need to know the multiplication maps 
$\mul_{21}$ and $\mul_{32}$, where we write $\mul_{mn}$ for the map
\begin{equation}
\label{equation4.0.5}
\mul_{mn} : \HoOC{m\D_0} \tensor \HoOC{n\D_0} \to \HoOC{(m+n) \D_0}.
\end{equation}
We then proceed as follows:
\begin{enumerate}
\item 
With respect to $\mul_{21}$, compute the space
$\HoOC{2\D_0 - \D}$ as 
\begin{equation}
\label{equation4.1}
\{ \s \in \HoOC{2\D_0} \mid \forall \t \in \HoOC{\D_0}, 
\quad \mul_{21} (\s \tensor \t) \in \WD \}.
\end{equation}
\item
Choose any nonzero $\f \in \HoOC{2\D_0 - \D}$, and note that $(f) = \D + \E
\linequiv 2\D_0$ (viewing $\f$ as a section of $\OC(2\D_0)$).  Then
using $\mul_{32}$, compute
\begin{equation}
\label{equation4.2}
\HoOC{5\D_0 - \D - \E} = \mul_{32} \left( \HoOC{3\D_0} \tensor \f \right) = 
\mul_{32} (\V \tensor \f).
\end{equation}
\item
Now use $\mul_{32}$ again to compute $\WE = \HoOC{3\D_0 - \E}$ as 
\begin{equation}
\label{equation4.3}
 \{ \s \in \V \mid
   \forall \t \in \HoOC{2D_0 - \D}, \quad \mul_{32}(\s \tensor \t) \in
			\HoOC{5\D_0 - \D - \E} \}.
\end{equation}
\end{enumerate}
\end{itemize}
\end{thmalg}

\begin{proof}
The first method proposed is completely straightforward.  As for justifying 
the second method, we note that the assertion in step~1 follows from the
usual argument, since our requirement that $\d_0 \geq 2\g$ implies that
$\HoOC{\D_0}$ has no base points.  Equation~\eqref{equation4.2} is
analogous to the identity in step~1 of Theorem/Algorithm~\ref{thmalg3.6}.
Finally, the assertion in step~3 follows since $\HoOC{2\D_0 - \D}$ also
has no base points.
\end{proof}

Our last major algorithm is to add two divisor classes on the Jacobian; by
now, the description is almost trivial.  We note a similar algorithm for
subtraction, which unexpectedly is slightly more efficient than the
addition algorithm.

\begin{thmalg}[Addition on the Jacobian]
\label{thmalg4.4}
Given small divisors $\D_1$ and $\D_2$, or rather their corresponding
subspaces $\W_{\D_1}$ and $\W_{\D_2}$, we can compute a space $\WE$
corresponding to a divisor $\E$ with $\xE = \x_{\D_1} + \x_{\D_2}$ (i.e.,
$\E \linequiv \D_1 + \D_2 - \D_0$) as follows:
\begin{enumerate}
\item
First use Proposition/Algorithm~\ref{propalg4.2.5} to compute the space
$\W_{\E'}$, where $\x_{\E'} = -(\x_{\D_1} + \x_{\D_2})$.
\item
Now use Theorem/Algorithm~\ref{thmalg4.3} to find a negation $\E$ of $\E'$.
\end{enumerate}
\end{thmalg}

\begin{propalg}[Subtraction]
\label{propalg4.5}
Given $\D_1$ and $\D_2$, we can find $\E$ such that $\x_{\D_1} - \x_{\D_2}
= \x_\E$ in {\em either one} of the following two ways:
\begin{itemize}
\item
First negate $\D_1$, thereby obtaining a small divisor $\D'$ with $x_{\D'}
= - \x_{\D_1}$.  Then apply Proposition/Algorithm~\ref{propalg4.2.5} to
$\D'$ and $\D_2$.  This yields $\WE$, where
$\xE = - ( -\x_{\D_1} + \x_{\D_2})$ as desired.  Note that this method only
involves one negation, instead of the two negations that would be involved
if we first negated $\D_2$, and then used the addition algorithm.
\item
Alternatively, we take the operations in Theorems/Algorithms
\ref{thmalg4.3} and~\ref{thmalg3.6.7}, and remove certain redundant steps
to streamline the subtraction algorithm.  This amounts to the following
steps:
\begin{enumerate}
\item
Compute $\HoOC{2\D_0 - \D_1}$ by ``dividing out'' $\D_0$, as in step~1 of
Theorem/Algorithm~\ref{thmalg4.3}.
\item
Take a nonzero $\f \in \HoOC{2\D_0 - \D_1}$, viewed as a section of
$\HoOC{2\D_0}$. Write $(\f) = \D_1 + \D'_1$, and use $\mul_{32}$ as in
step~2 of Theorem/Algorithm~\ref{thmalg4.3} to calculate
\begin{equation}
\label{equation4.4}
\HoOC{5\D_0 - \D_1 - \D'_1 - \D_2} =
% \mul_{32} \left( \HoOC{3\D_0 - \D_2} \tensor \f \right) = 
\mul_{32} (\W_{\D_2} \tensor \f).
\end{equation}
\item
Analogously to step 3 in  Theorem/Algorithm~\ref{thmalg4.3}, compute
$\HoLm{\D'_1 - \D_2} = \HoOC{3\D_0 - \D'_1 - \D_2}$ as
\begin{equation}
\label{equation4.5}
 \{ \s \in \V \mid
   \forall \t \in \HoOC{2D_0 - \D_1}, \quad \mul_{32}(\s \tensor \t) \in
			\HoOC{5\D_0 - \D_1 - \D'_1 - \D_2} \}.
\end{equation}
\item
Finally, apply Theorem/Algorithm~\ref{thmalg3.6} to obtain $\WE$ as the
``flip'' of the space $\HoLm{\D'_1 - \D_2}$.
\end{enumerate}
\end{itemize}
\end{propalg}
\begin{proof}
The first method is completely straightforward.  As for the second, it
boils down to the facts that $\x_{\D'_1} = -\x_{\D_1}$ and that $\xE = -(
\x_{\D'_1} + \x_{\D_2})$.
\end{proof}

\section{Further improvements to the algorithms in Section~\ref{section4}:
the medium and small models}
\label{section5}

We now sketch two ways to modify the algorithms from the previous section,
while speeding up the algorithms by a constant factor, so the asymptotics
of our algorithms become a smaller multiple of $O(\g^4)$.  In both cases,
we achieve this linear speedup by reducing the degree of the basic line
bundle $\L$, and with it the dimension of the ambient space $\V$.

We shall refer to our first method as the {\em medium model}.  We fix as
before an effective divisor $\D_0$ of degree $\d_0 \geq 2\g + 1$.  (It is
possible, but a bit awkward, to work with $\d_0 = 2\g$, but this would
hamper us in our liberal use of Lemma~\ref{lemma2.1}.)  We let $\L =
\OC(2\D_0)$, and represent a point $\xD = [\D - \D_0]$ by the space $\WD$,
as before.  Here $\deg \D = \d_0$; in case $\d_0 = 2\g + 1$, we obtain that
$\WD$ is a $(\g+2)$-dimensional subspace of the $(3\g + 3)$-dimensional
space $\V = \HoOC{2\D_0}$.  Since $\thedim = \dim \V = 3\g + 3$ in this
case, we have improved the quantity $\thedim^4$ by a factor of around
$(5/3)^4 \approx 7.7$ over the previous section; this is offset by the fact
that we now have to do about twice as much linear algebra per operation in
the Jacobian.  So we predict an effective speedup by approximately a factor
of $4$.

Our first observation regarding the medium model is that we can still use
Theorems/Algorithms \ref{thmalg3.7} and~\ref{thmalg4.1} to test for
membership of and equality in the Jacobian in the medium model.  We also
observe that negation is now simply a matter of ``flipping,'' as in
Theorem/Algorithm~\ref{thmalg3.6}.  We need to make some modifications to
the addition and subtraction algorithms, though.  We first give the
algorithm for the ``addflip'' operation, analogously to
Proposition/Algorithm~\ref{propalg4.2.5}.  The idea is virtually the same,
but we need to keep track of more multiplication maps.  We shall generally
use the notation $\mul_{mn}$ to denote the same multiplication map as in
equation~\eqref{equation4.0.5}:
\begin{equation}
\label{equation5.1}
\mul_{mn} : H^0(m\D_0) \tensor H^0(n\D_0) \to H^0( (m+n) \D_0 ).
\end{equation}
Note that we have simplified notation in this section by writing $m\D_0$
instead of, say, $\OC(m\D_0)$.

\begin{propalg}[Addflip (medium model)]
\label{propalg5.1}
Given $\W_{\D_1} = H^0(2\D_0 - \D_1)$ and $\W_{\D_2} = H^0(2\D_0
- \D_2)$, we calculate $\WE$, where $\xE = - (\x_{\D_1} + \x_{\D_2})$, as
follows:
\begin{enumerate}
\item
Compute $H^0(4\D_0 - \D_1 - \D_2) = \mul_{22}( \W_{\D_1} \tensor
\W_{\D_2})$, by Lemma~\ref{lemma2.1}.
\item
Compute $H^0(3\D_0 - \D_1 - \D_2) = \{ \s \in H^0(3\D_0) \mid
	\forall \t \in H^0(\D_0), \quad
	\mul_{31} (\s \tensor \t) \in H^0(4\D_0 - \D_1 - \D_2) \}$.
\item
Take a nonzero $\f \in H^0(3\D_0 - \D_1 - \D_2)$; viewing $\f$ as a section 
of $\OC(3\D_0)$, we obtain $(\f) = \D_1 + \D_2 + \E$, where $\xE$ is the
desired class.  Now compute
$H^0(5\D_0 - \D_1 - \D_2 - \E) = \mul_{32}(\f \tensor H^0(2\D_0))$.
\item
Obtain $\WE = \{ \s \in \V \mid
	\forall \t \in H^0(3\D_0 - \D_1 - \D_2), \quad
	\mul_{23}(\s \tensor \t) \in H^0(5\D_0 - \D_1 - \D_2 - \E) \}$.
\end{enumerate}
\end{propalg}

The addition algorithm is then obtained by performing an addflip, followed
by a regular flip to negate the result.  Subtraction can be done similarly
to the first method of Proposition/Algorithm~\ref{propalg4.5}.
Alternatively, one can adapt the second method from
Proposition/Algorithm~\ref{propalg4.5}: step~1 in that setting is
superfluous, and one can obtain $H^0(3\D_0 - \D_2) = \mul_{21}(\W_{\D_2}
\tensor H^0(\D_0))$ for use in step~2 there, and continue with step~3.
Then modify step~4 of that algorithm to resemble the last two steps of
Proposition/Algorithm~\ref{propalg5.1}.  We have described the alternative
subtraction algorithm to illustrate possible techniques, but the simpler
subtraction is probably just as quick.

We next turn to what we shall call the {\em small model}.  In principle, one 
can represent all points on the Jacobian using divisors of degree $\g$,
which can be uniquely described using an ambient line bundle of degree
$3\g$.  However, in order to use Lemma~\ref{lemma2.1},
and to avoid some contortions, we shall use instead a divisor $\D_0$ of
degree $\d_0 = \g + 1$ (or at least $\d_0 \geq \g + 1$), and put $\L =
\OC(3\D_0)$.  So in this case $\V = \HoOC{3\D_0}$ has dimension $2\g + 4$,
but we need to go up to at least $H^0(7\D_0)$ within our calculation.  We
use the same notation for $\mul$ as in~\eqref{equation5.1}.  The algorithms
in this setting are reasonably direct generalizations of those in
Section~\ref{section4} and of those in the medium model; the main problem
is that a divisor $\D'$ of degree $2\d_0$ may not be recognizable from the
space $\HoLm{\D'}$, which may have base points; so we need to work a bit
more to deal with $\D'$ by looking instead at $H^0(\L'-\D')$ for a line
bundle $\L'$ of sufficiently high degree.  On the other hand, the small
model uses smaller-dimensional vector spaces and should give the fastest
algorithms among the ones discussed in this paper.  In the small model, we
can test equality on the Jacobian using the second method of
Theorem/Algorithm~\ref{thmalg4.1} (the reader is invited to find an analog
of the first method, using techniques analogous to those in
Propositions/Algorithms \ref{propalg5.2}--\ref{propalg5.4}).  We
therefore content ourselves with the algorithms to test for membership of
the Jacobian, for the addflip operation, and a slight improvement for
negation (instead of reducing it to addflip).  Addition and subtraction
can then be done in a straightforward way, as in Section~\ref{section4}.
We note that one can obtain a slight improvement to the subtraction
algorithm, similarly to the second method in
Proposition/Algorithm~\ref{propalg4.5}.  The details are slightly lengthy,
and are left to the reader.

\begin{propalg}[Membership test (small model)]
\label{propalg5.2}

Let $\W \subset \V$ be a codimension $\d_0$ subspace.  Then we can test
whether $\W = \WD$ for some $\D$ as follows:
\begin{enumerate}
\item 
Let $\f$ be a nonzero element of $\W$, and form the space $\W_+ = \mul_{43} 
(H^0(4\D_0) \tensor \f) \subset H^0(7\D_0)$.
\item
Calculate $\W' = \{ \s \in H^0(4\D_0) \mid \forall \t \in \W, \quad
\mul_{43}(\s \tensor \t) \in \W_+ \}$.
\item
Then $\W$ comes from a divisor $\D$ if and only if $\W'$ has codimension
$2\d_0$ in $H^0(4\D_0)$.
\end{enumerate}
\end{propalg}
\begin{proof}
This is entirely analogous to Theorem/Algorithm~\ref{thmalg3.7}.
The idea is that we hope to have $(\f) = \D + \D'$, whence $\W_+ =
H^0(7\D_0 - \D - \D')$ and $\W' = H^0(4\D_0 - \D')$.  We have had to use
the higher-degree line bundle $\OC(4\D_0)$ in addition to $\L$, in
order to ensure that the degrees are large enough for us to be able to use
Lemma/Algorithm~\ref{lemmalg2.2}.
\end{proof}

\begin{propalg}[Addflip (small model)]
\label{propalg5.3}
Given two subspaces $\W_{\D_1} = H^0(3\D_0 - \D_1)$ and $\W_{\D_2} =
H^0(3\D_0 - \D_2)$, we can find $\WE$, where $\xE = -(\x_{\D_1} +
\x_{\D_2})$, by modifying Proposition/Algorithm~\ref{propalg4.2.5} as
follows: 
\begin{enumerate}
\item
Compute $H^0(3\D_0 - \D_1 - \D_2)$ by Theorem/Algorithm~\ref{thmalg3.4},
without discarding the intermediate result $H^0(6\D_0 - \D_1 - \D_2)$.
(It is a good idea to obtain $H^0(4\D_0 - \D_1 - \D_2)$ in the process;
this can be done by setting up linear equations to simultaneously
``divide'' out both $H^0(3\D_0)$ and $H^0(2\D_0)$ from $H^0(6\D_0 - \D_1
- \D_2)$, the former via \eqref{equation3.3}, and the latter via an
analogous calculation that largely overlaps with \eqref{equation3.3}.)
\item
Take a nonzero $\f \in H^0(3\D_0 - \D_1 - \D_2)$; so, viewing $\f$ as a
section of $\OC(3\D_0)$, we have $(\f) = \D_1 + \D_2 + \E$ for the desired
$\E$. 
\item
Calculate $\WE$ as
\begin{equation}
\label{equation5.2}
 \{ \s \in H^0(3\D_0) \mid \forall \t \in H^0(6\D_0 - \D_1 - \D_2), \quad 
\mul_{36}(\s \tensor \t) \in \mul_{36} (\f \tensor H^0(6\D_0)) \}.
\end{equation}
Note that this method involves calculation in $H^0(9\D_0) = H^0(3\L)$; if
we have done the extra computation in step~1, we can limit ourselves to
$H^0(7\D_0)$, by replacing $6\D_0$ by $4\D_0$ throughout
\eqref{equation5.2}.
\end{enumerate}
\end{propalg}
\begin{proof}
Immediate.  Note the trouble that we went to because $H^0(3\D_0 -
\D_1 - \D_2)$ is not necessarily base point free, which meant that we had
to use \eqref{equation5.2} instead of doing a regular ``flipping''
algorithm after step~1 above.  However, we still point out that,
generically, $H^0(3\D_0 - \D_1 - \D_2)$ {\em is} base point free, since its
degree is $\d_0 \geq \g + 1$.
\end{proof}

\begin{propalg}[Negation (small model)]
\label{propalg5.4}
Given $\WD$, we can find $\WE$, where $\xE = -\xD$, by using the addflip
operation.  Alternatively, we can use the following slightly faster
method:
\begin{enumerate}
\item
Calculate $H^0(5\D_0 - \D) = \mul_{23}( H^0(2\D_0) \tensor \WD )$.  Then
calculate $H^0(2\D_0 - \D) = \{ \s \in H^0(2\D_0) \mid \forall \t \in \V,
\quad \mul_{23} (\s \tensor \t) \in H^0(5\D_0 - \D) \}$.
\item
Take a nonzero $\f \in H^0(2\D_0 - \D)$, and compute $H^0(6\D_0 - \D - \E)
= \mul_{24}( \f \tensor H^0(4\D_0))$.  
\item 
Finally, compute $\WE = \{ \s \in H^0(3\D_0) \mid \forall \t \in \WD,
\quad \mul_{33} (\s \tensor \t) \in H^0(6\D_0 - \D - \E) \}$.
\end{enumerate}
\end{propalg}
\begin{proof}
This is just like Theorem/Algorithm~\ref{thmalg4.3}, but we have included
the algorithm to show in step~1 the technique of going from $H^0(3\D_0 -
\D)$ to $H^0(2\D_0 - \D)$ by virtue of first ``going up,'' then ``going
down.''  This is because the degree of $\D_0$ is too small for us to
automatically be able to ``divide'' by $H^0(\D_0)$.  If we were originally
careful enough to select $\D_0$ so that $H^0(\D_0)$ were base point free,
we could have begun as in Theorem/Algorithm~\ref{thmalg4.3}.  However, we
cannot in any way guarantee that $H^0(2\D_0 - \D)$ is base point
free, since $\D$ is arbitrary; hence, concluding with the computation
in~\eqref{equation4.3} was never an option in the setting of the small
model.
\end{proof}

We conclude this paper with a brief discussion of converting points on the
Jacobian between their representations in the different models.  The first
observation is that it is easy to pass between two representations of the
same (as always, effective) divisor $\D$ in terms of spaces $\WD$ of
sections of two different line bundles $\L_1$ and $\L_2$.  (Typically, $\L_1
= \OC(n_1 \D_0)$ and $\L_2 = \OC(n_2 \D_0)$ for a fixed $\D_0$ and integers
$n_1, n_2$.)  So the problem is to pass from a knowledge of $H^0(\L_1 - \D)$
to one of $H^0(\L_2 - \D)$.  Let us assume for a moment that $\deg \L_1
\geq 2\g + 1 + \deg \D$.  Then take auxiliary line bundles $\L'_1, \L'_2$
of degree at least $2\g + 1$, such that $\L_1 + \L'_1 = \L_2 + \L'_2$;
here, addition of line bundles really refers to their tensor product.
Then, using multiplication of sections in $H^0(\L_1 - \D)$ and
$H^0(\L'_1)$, we can obtain the space $H^0(\L_1 + \L'_1 - \D) = H^0(\L_2 +
\L'_2 - \D)$.  We can then ``divide'' by $H^0(\L'_2)$ to obtain $H^0(\L_2 -
\D)$, as desired.  This is similar to step~1 in
Proposition/Algorithm~\ref{propalg5.4}.  Of course, if the degree of $\L_1$ 
is sufficiently larger than the degree of $\L_2$, we only need to perform a
single division, without first raising the degree using $\L'_1$.
Similarly, if the degree of $\L_2$ is large compared to that of $\L_1$,
then a single multiplication step is sufficient.

The above discussion seems to require the use of Lemma~\ref{lemma2.1} for
the multiplication step.  We can modify the procedure, though, to work
whenever $H^0(\L_1 - \D)$ is base point free (such as when $\deg \L_1 \geq
2\g + \deg \D$).  The idea is to take a nonzero section $\f \in H^0(\L_1 -
\D)$, whose divisor (viewing $\f$ as a section of $\L_1$) is $\D + \D'$.
Here we caution that the degree of $\D'$ may be large, so that $H^0(\L_1 -
\D')$ may well have base points.  To get around this, we take an auxiliary
line bundle $\L'$ of degree at least $2\g + 1 + \deg \D'$, and multiply
$\f$ by $H^0(\L')$ to obtain $H^0(\L_1 + \L' - \D - \D')$.  Dividing out
$H^0(\L_1 - \D)$ then yields $H^0(\L' - \D')$.  We can then multiply $\f$
by $H^0(\L'_1 + \L')$, where $\L'_1$ has suitably large degree, and then
divide by $H^0(\L' - \D')$ (which is now base point free), to obtain
$H^0(\L'_1 + \L_1 - \D)$ even in the case where we could not use
Lemma~\ref{lemma2.1} as in the above paragraph.  We note that throughout
all these computations, the degrees of the intermediate line bundles are at
most $O(\max(\deg \L_1, \deg \L_2, \g))$.

Now that we have full control over the ambient line bundle $\L$ used to 
represent individual divisors $\D$, we shall always assume that the degree
of $\L$ is sufficiently large, and use $\L$ to represent all divisors
(so $\D$ is given by $\WD = H^0(\L - \D)$, unless stated otherwise).  We
now discuss how to change the basepoint $\D_0$ used to map a divisor $\D$
to the point $\xD = [\D - \D_0]$ on the Jacobian.  It is always easy to add
an effective divisor $\E$ to $\D_0$, since $[\D - \D_0] = [(\D + \E) -
(\D_0 + \E)]$, so we can use Theorem/Algorithm~\ref{thmalg3.4} to
replace $\WD$ with $\W_{\D + \E}$.  In particular, we can also assume that
$\deg \D_0 \geq 2\g + 1$.  Thus the main question that remains is how to
replace $\D_0$ by a new basepoint $\E_0$ of smaller degree.  This means
that given $\D$, we wish to find $\E$ such that $[\D - \D_0] = [\E -
\E_0]$; in other words, we seek an effective divisor $\E \linequiv \D +
\E_0 - \D_0$.  Writing $\D_1 = \D + \E_0$ (which we compute as usual by
Theorem/Algorithm~\ref{thmalg3.4}, and whose degree is even larger
than $\deg \D = \deg \D_0 \geq 2\g + 1$), we have reduced our question to
the following result, which leads to an efficient algorithm for the
Riemann-Roch problem on $\C$.  Like our other algorithms, this one requires 
$O(\max(\deg \D_0, \deg \D_1, \g)^4)$ field operations in $\k$.  Remember
that there is no loss of generality in letting $\L$ have large degree $\N
\geq 2\g + 1 + \deg \D_1$.

\begin{thmalg}
\label{thmalg5.5}
Let $\deg \L = \N$, and assume given (effective) divisors $\D_1$ and
$\D_0$, whose degrees lie between $2\g + 1$ and $\N - 2\g - 1$.  Then we
can compute the space $H^0(\D_1 - \D_0)$, and with it (if the space is
nonzero) an effective divisor $\E \linequiv \D_1 - \D_0$.  The procedure
is:
\begin{enumerate}
\item
Use Theorem/Algorithm~\ref{thmalg3.6} to compute $\D'_1$ (specifically, to
compute $\W_{\D'_1} = H^0(\L - \D'_1)$), where $\D_1 + \D'_1 \linequiv \L$.
Thus $\L - \D'_1 - \D_0 \linequiv \D_1 - \D_0$, and an explicit isomorphism
from $H^0(\D_1 - \D_0)$ to $H^0(\L - \D'_1 - \D_0)$ is given by
multiplication by the section $\f$ of $\L$ with 
$(\f) = \D_1 + \D'_1$ used in the calculation of $\W_{\D'_1}$.
\item
Compute $H^0(2\L - \D'_1 - \D_0) = \mul( \W_{\D'_1} \tensor \W_{\D_0})$,
and use it to compute $\W_{\D'_1 + \D_0} \isomorphic H^0(\D_1 - \D_0)$, as in
Theorem/Algorithm~\ref{thmalg3.4}.
\item
A divisor $\E$ as above exists if and only if we can find a nonzero $\g \in 
\W_{\D'_1 + \D_0}$, in which case we view $\g$ as a section of $\L$ and see 
that its divisor is $(\g) = \D'_1 + \D_0 + \E$.  Then calculate
\begin{equation}
\label{equation5.3}
\WE = \{ \s \in \HoL \mid \forall \t \in H^0(2\L - \D'_1 - \D_0), \quad
\mul'(\s \tensor \t) \in \mul'(\g \tensor H^0(2\L)) \}.
\end{equation} 
Here $\mul': H^0(\L) \tensor H^0(2\L) \to H^0(3\L)$ is the multiplication
map.  (We did not use the notational convention of \eqref{equation5.1},
since we are not assuming that $\L$ is a multiple of $\D_0$.)
\end{enumerate}
\end{thmalg}
\begin{proof}
We know that $\mul'(\g \tensor H^0(2\L)) = H^0(3\L - \D'_1 - \D_0 - \E)$,
and our generous assumptions about degrees ensure that $H^0(2\L - \D'_1 -
\D_0)$ is base point free.  Note that we do not assume that the smaller
space $H^0(\L - \D'_1 - \D_0)$ is base point free; indeed, its degree is
$\deg \D_1 - \deg \D_0$, which may be quite small.  The rest of the proof
is standard by now.
\end{proof}

\bibliographystyle{amsalpha}

\end{document}